\documentclass[11pt,reqno]{amsart}
\usepackage{amsmath}
\usepackage{amssymb}
\usepackage{amstext}
\usepackage{mathrsfs}
\usepackage{a4wide}
\usepackage{graphicx}
\usepackage{amsfonts}
\usepackage{mathtools}
\usepackage{esint}
\usepackage{color}
\usepackage{hyperref}
\usepackage{enumerate}
\usepackage{stmaryrd}

\allowdisplaybreaks \numberwithin{equation}{section}
\hypersetup{colorlinks=true, linkcolor=blue, filecolor=blue, urlcolor=blue, citecolor=red}

\newtheorem{theorem}{Theorem}[section]

\newtheorem{lemma}[theorem]{Lemma}
\theoremstyle{definition}

\theoremstyle{remark}
\newtheorem{remark}[theorem]{Remark}

\newcommand{\ds}{\displaystyle}
\newcommand{\R}{\mathbb{R}}
\newcommand{\C}{\mathbb{C}}
\newcommand{\N}{\mathbb{N}}
\newcommand{\T}{\mathbb{T}}
\newcommand{\Z}{\mathbb{Z}}
\newcommand{\ep}{\varepsilon}

\newcommand{\Pm}{\mathcal{P}}
\newcommand{\dd}{\,\mathrm{d}}
\newcommand{\dy}{\,\mathrm{d}y}
\newcommand{\dx}{\,\mathrm{d}x}
\newcommand{\dt}{\,\mathrm{d}t}
\newcommand{\dtau}{\,\mathrm{d}\tau}

\newcommand{\bfe}{\boldsymbol{e}}
\newcommand{\vv}{\mathbf{v}}
\newcommand{\bfF}{\mathcal{F}}
\newcommand{\bfG}{\mathcal{G}}
\newcommand{\gfr}{g}
\newcommand{\hfr}{h}

\newcommand{\norm}[1]{\left\Vert#1\right\Vert}

\newcommand{\lot}{\mathrm{l.o.t.}}
\newcommand{\Dself}[1]{\varepsilon^2b_{#1}^2\big(r_{#1}(x)-r_{#1}(y)\big)^2+4\big(1+\varepsilon b_{#1}r_{#1}(x)\big)\big(1+\varepsilon b_{#1}r_{#1}(y)\big)\sin^2\big(\tfrac{x-y}{2}\big)}
\newcommand{\Dinter}[1]{\big(\varepsilon b_{#1}R_{#1}(x)\cos x+\varepsilon b_{3-#1}R_{3-#1}(y)\cos y-d\big)^2+\big(\varepsilon b_{#1}R_{#1}(x)\sin x+\varepsilon b_{3-#1}R_{3-#1}(y)\sin y\big)^2}

\newcommand{\Ufull}[1]{\big(1+\varepsilon b_{#1}r_{#1}(x)\big)^2+\varepsilon^2b_{#1}^2r_{#1}'(x)^2}
\newcommand{\DhA}[2]{2\varepsilon^2b_{#1}^2\big(r_{#1}(x)-r_{#1}(y)\big)\big(#2(x)-#2(y)\big)}
\newcommand{\DhB}[2]{4\varepsilon b_{#1}\Big(#2(x)\big(1+\varepsilon b_{#1}r_{#1}(y)\big)+\big(1+\varepsilon b_{#1}r_{#1}(x)\big)#2(y)\Big)\sin^2\big(\tfrac{x-y}{2}\big)}

\newcommand{\Dq}[1]{\frac{r_{#1}(x)-r_{#1}(y)}{2\sin(\frac{x-y}{2})}}
\newcommand{\Dn}[1]{\big(1+\varepsilon b_{#1}r_{#1}(x)\big)\big(1+\varepsilon b_{#1}r_{#1}(y)\big)+\varepsilon^2b_{#1}^2\Big(\Dq{#1}\Big)^2}
\newsavebox{\fitbox}
\newcommand{\fit}[2][\hsize]{\sbox{\fitbox}{$\displaystyle #2$}\ifdim\wd\fitbox>#1\resizebox{#1}{!}{\usebox{\fitbox}}\else\usebox{\fitbox}\fi}

\begin{document}

\title[Asymmetric vortex sheets]{Co-rotating and traveling asymmetric vortex sheets for the planar Euler equation}

\author{Edison Cuba}
\address{Departamento Académico de Matemáticas y Estadística, Universidad Nacional de San
Antonio Abad del Cusco, Perú.}
 \email{edison.cuba@unsaac.edu.pe}


\begin{abstract}
The aim of this paper is to construct co-rotating and traveling vortex sheet pairs with asymmetric configurations for the two-dimensional incompressible Euler equation. The vorticity of these solutions is concentrated on two small closed curves of different sizes, which carry vortex sheets of different circulations and either rotate uniformly about their centroid or translate with constant speed. The proof relies on the Birkhoff--Rott formulation of the problem and on the implicit function theorem at the point vortex pairs.
\end{abstract}

\maketitle

% \tableofcontents 

\section{Introduction}

\subsection{Vortex sheets for the two-dimensional Euler equation}

The vorticity formulation of the two-dimensional incompressible Euler equation in the whole plane reads
\begin{equation}\label{1-1}
\begin{cases}
\partial_t\vartheta+\vv\cdot\nabla\vartheta=0 & \text{in}\ \R^2\times(0,T),\\
\vartheta\big|_{t=0}=\vartheta_0 & \text{in}\ \R^2,
\end{cases}
\end{equation}
where $\vv=(v_1,v_2)$ is the velocity field and $\vartheta=\partial_1v_2-\partial_2v_1$ is the vorticity. By the incompressibility condition $\nabla\cdot\vv=0$, the velocity is recovered from the vorticity through the Biot--Savart law
\begin{equation*}
\vv=\nabla^\perp(-\Delta)^{-1}\vartheta,\qquad (-\Delta)^{-1}\vartheta(x)=\frac{1}{2\pi}\int_{\R^2}\ln\frac{1}{|x-y|}\,\vartheta(y)\dd y,
\end{equation*}
where,  $(x_1,x_2)^\perp:=(x_2,-x_1)$. Equivalently,
\begin{equation*}
\vv(x)=-\frac{1}{2\pi}\int_{\R^2}\frac{(x-y)^\perp}{|x-y|^2}\,\vartheta(y)\dd y .
\end{equation*}
The equation \eqref{1-1} expresses the fact that the vorticity is transported by the flow. For bounded and integrable initial vorticity, the global existence and uniqueness of weak solutions is due to Yudovich \cite{Yud}. This class contains the vortex patches, namely the characteristic functions of bounded domains. See \cite{MB,MP} for more details.

A \emph{vortex sheet} is a weak solution of \eqref{1-1} whose vorticity is a measure concentrated on curves. Physically, vortex sheets model a surface of discontinuity of the tangential velocity, which is generated when a fluid of small viscosity passes a sharp edge or a rigid wall, or at the interface between two streams with different velocities. In these situations the viscosity plays an important role in the generation of the sheet, where its evolution is well described by the Euler equation \cite{Bir,MB,Saf}. In this paper the vorticity $\vartheta(t)$ will be supported on two simple closed curves $\Gamma_i(t)=\{z_i(t,x):\ x\in[0,2\pi)\}$, $i=1,2$, namely
\begin{equation}\label{measure}
\int_{\R^2}\varphi\dd\vartheta(t)=\sum_{i=1}^{2}\int_0^{2\pi}\varphi(z_i(t,x))\,\omega_i(t,x)\dx,\qquad \forall\,\varphi\in C^\infty_c(\R^2),
\end{equation}
where,  $\omega_i(t,x)$ is the vortex strength with respect to the parametrization $z_i$, $\delta_{\Gamma}$ is the Dirac measure on the curve $\Gamma$, and $\int_0^{2\pi}\omega_i\dx$ is the circulation of $\Gamma_i$ and
\begin{equation*}
\gamma_i(t,x):=\frac{\omega_i(t,x)}{|\partial_xz_i(t,x)|}
\end{equation*}
is the regularized vortex strength of $\Gamma_i$. That is, the strength with respect to the arc length.  The velocity field is smooth away from the curves and has a jump across each $\Gamma_i$, its normal component is continuous, while the jump of its tangential component is equal to the strength with respect to the arc length. On the sheet we use the average of the two one-sided limits, which is given by the Birkhoff--Rott integral
\begin{equation}\label{1-2}
\begin{split}
\vv(t,z_i(t,x))&=-\sum_{j=1}^2\mathcal{BR}(z_j,\omega_j)(z_i(t,x))\\
&:=-\frac{1}{2\pi}\sum_{j=1}^{2}PV\!\!\int_0^{2\pi}\frac{(z_i(t,x)-z_j(t,y))^\perp}{|z_i(t,x)-z_j(t,y)|^2}\,\omega_j(t,y)\dy,
\end{split}
\end{equation}
where the principal value is only needed for the case $j=i$. The evolution of the sheets is then described by the Birkhoff--Rott system (see \cite{Bir,CQZ,MB,Rott,Saf})
\begin{align}
\partial_t z_i(t,x)&=\vv(t,z_i(t,x))+c_i(t,x)\,\partial_xz_i(t,x),\label{1-3}\\
\partial_t\omega_i(t,x)&=\partial_x\big(c_i(t,x)\omega_i(t,x)\big),\label{1-4}
\end{align}
where $c_i$ is the reparametrization freedom. Since only the curve and the measure \eqref{measure} have a physical meaning, the tangential component of the velocity of the parametrization can be chosen arbitrarily, and \eqref{1-4} expresses the conservation of the circulation along the curve. For simplicity we omit the notation $PV$ in what follows.

The mathematical analysis of \eqref{1-3}--\eqref{1-4} is highly challenging. A constant-strength flat sheet is a steady state but is subject to Kelvin–Helmholtz instability, making the Cauchy problem ill-posed in Sobolev spaces \cite{CO, Ebin}. However, local existence holds for analytic data \cite{SS, SSBF}.
Global solutions exist for small analytic perturbations \cite{Duc}, and Lebeau \cite{Leb} analyzed their regularity. Conversely, finite-time curvature blow-up can occur even for initially analytic sheets \cite{Kra, Moo}. Wu \cite{Wu} established local existence for chord-arc curves with weakly regular strength. Finally, global weak solutions exist for single-signed vorticity \cite{Del, Maj} and under reflection symmetry \cite{LNX}.

\subsection{Relative equilibria and desingularization of point vortices}

Motivated by the instability of these configurations, identifying nontrivial relative equilibria. Specifically, rigidly translating or rotating vortex sheets constitutes a non-trivial problem.  The classic point vortex system serves as an excellent source for such states. Assuming the vorticity is concentrated as a finite sum of Dirac delta functions, $\sum_j2\pi\gamma_j\boldsymbol{\delta}_{Z_j(t)}$, we can formally reduce the Euler equations to the well-known Kirchhoff--Helmholtz system:
\begin{equation*}
\frac{\dd Z_j}{\dd t}=-\sum_{k\neq j}\gamma_k\frac{(Z_j-Z_k)^\perp}{|Z_j-Z_k|^2},
\end{equation*}
see \cite{Aref,MP,New,Saf}. For two point vortices of circulations $2\pi\gamma_1$ and $2\pi\gamma_2$ placed at $0$ and $d\bfe_1$, $d>0$, one easily checks that the pair rotates uniformly about its centroid $(x_0^*,0)$, provided $\gamma_1+\gamma_2\neq0$, with
\begin{equation}\label{ang}
\Omega^*:=\frac{\gamma_1+\gamma_2}{d^2},\qquad x_0^*:=\frac{d\gamma_2}{\gamma_1+\gamma_2},
\end{equation}
while, if $\gamma_2=\gamma_1$, the point vortices of circulations $2\pi\gamma_1$ and $-2\pi\gamma_2$ translate in the direction $\bfe_2$ with speed
\begin{equation}\label{tra}
U^*:=\frac{\gamma_1}{d}.
\end{equation}

Relative equilibria for vortex patches have been extensively studied. Rotating patches (V-states) were discovered numerically \cite{DZ} and proven rigorously near the disk via bifurcation theory \cite{Bur}, with later extensions addressing regularity, doubly connected domains, and global bifurcation \cite{GHS,HMW,HMV,dHHMV}.  Co-rotating and translating pairs of equal patches were computed in \cite{Pie,SSz} and established variationally in \cite{CLW,SVS,Tur}. Alternatively, Hmidi and Mateu \cite{HM} developed a contour dynamics approach that was subsequently adapted for the generalized surface quasi-geostrophic (gSQG) equations \cite{CQSZ,CCG2,CCG1}. Physically relevant asymmetric configurations of varying sizes and circulations, initially explored numerically \cite{Dri}, have also been rigorously constructed. Hassainia and Hmidi \cite{HH} proved the existence of co-rotating and traveling asymmetric pairs for the Euler equations, while analogous multipole and asymmetric patch equilibria for gSQG were established in the ranges $\alpha\in[0,1)$ by \cite{HW} and $\alpha\in[1,2)$ by \cite{CF}.

Our understanding of relative equilibria for vortex sheets,  remains much more restricted. Apart from circles with uniform strength and segments, Protas and Sakajo \cite{PS} found rotating equilibria made of segments, and G\'omez-Serrano, Park, Shi and Yao obtained uniformly rotating sheets bifurcating from the circle \cite{GPSY2}, together with rigidity results showing that stationary or uniformly rotating sheets with positive strength must be radial in some regimes \cite{GPSY1,GPSY3}. In \cite{CQZ}, the authors constructed a new class of co-rotating and traveling vortex sheets supported on several small closed curves, by applying the implicit function theorem to the Birkhoff--Rott formulation at point vortex configurations. 

\subsection{Main results}

The purpose of this paper is to extend the results of \cite{CQZ} to the asymmetric configurations. Let $d>0$, $b_1,b_2>0$ and $\kappa_1,\kappa_2\in\R\setminus\{0\}$. We look for vortex sheets supported on two closed curves close to the circles of radii $\ep b_1$ and $\ep b_2$ centered at the origin and at $d\bfe_1$, respectively, with circulations $2\pi\kappa_1$ and $2\pi\kappa_2$, which either rotate uniformly about their centroid $(x_0,0)$ or travel with constant speed $U$ in the vertical direction. The initial vorticity is of the form
\begin{equation*}
\vartheta_{0,\ep}=\frac{\kappa_1}{\ep b_1}\omega_1\delta_{\Gamma^\ep_1}\pm\frac{\kappa_2}{\ep b_2}\omega_2\delta_{\Gamma^\ep_2},
\end{equation*}
where,  $\omega_i$ is the vortex strength of the curve $\Gamma_i^\ep$ with respect to its parametrization $z_i$, and the factor $\frac{\kappa_i}{\ep b_i}$ compensates the size of the curve. As $\ep\to0$, these configurations converge to the pairs of point vortices described in \eqref{ang} and \eqref{tra}. The following is an informal version of our main result. See Theorems \ref{main} and \ref{trans} for the precise statements.

\begin{theorem}\label{thm:informal}
Let $d>0$ and $b_1,b_2>0$.
\begin{enumerate}
\item[(1)] For any $\kappa_1,\kappa_2\in\R\setminus\{0\}$ with $\kappa_1+\kappa_2\neq0$, there exists $\ep_0>0$ such that for every $\ep\in(0,\ep_0]$, \eqref{1-1} has a co-rotating asymmetric vortex sheet pair
$\vartheta_\ep(t,X)=\vartheta_{0,\ep}\big(Q_{-\Omega t}(X-x_0\bfe_1)+x_0\bfe_1\big)$, where $\vartheta_{0,\ep}$ is given in \eqref{sole}, the two sheets have circulations $2\pi\kappa_1+O(\ep^4)$ and $2\pi\kappa_2+O(\ep^4)$, and
\begin{equation*}
\Omega=\Omega^*+O(\ep),\qquad x_0=x_0^*+O(\ep).
\end{equation*}
\item[(2)] For any $\kappa_1\in\R\setminus\{0\}$, there exist $\ep_0>0$ and $\kappa_2=\kappa_2(\ep)=\kappa_1+O(\ep)$ such that for every $\ep\in(0,\ep_0]$, \eqref{1-1} has a traveling asymmetric vortex sheet pair $\vartheta_\ep(t,X)=\vartheta_{0,\ep}(X-tU\bfe_2)$, where $\vartheta_{0,\ep}$ is given in \eqref{initial}, the two sheets have circulations $2\pi\kappa_1+O(\ep^4)$ and $-2\pi\kappa_2(\ep)+O(\ep^4)$, and $U=U^*+O(\ep)$.
\end{enumerate}
In both cases, each component of the solution is supported on a $C^\infty$ closed curve with
convex interior.
\end{theorem}

\begin{remark}\label{rem:intro}
\begin{enumerate}
\item[(a)] When $\kappa_1=\kappa_2$ and $b_1=b_2$, the uniqueness part of Theorem \ref{main} shows that $x_0=d/2$ and that the two sheets are images of each other under the rotation of angle $\pi$ about $(d/2,0)$, so we recover the case $m=2$ of \cite[Theorem 1.1]{CQZ} (where the distance between the centers is denoted by $2d$). Similarly, when $b_1=b_2$ in Theorem \ref{trans} one has $\kappa_2(\ep)=\kappa_1$, and we recover \cite[Theorem 1.2]{CQZ}.
\item[(b)] The condition $\kappa_1+\kappa_2\neq0$ in (1) is natural. It is exactly the condition under which the pair of point vortices rotates about a finite centroid, and it is also the condition under which the linearized operator is an isomorphism, see Lemma \ref{lem:iso}. 
\end{enumerate}
\end{remark}

\medskip

The paper is organized as follows. In Section \ref{sec2} we introduce the functional spaces, and we derive the equations governed by the co-rotating and traveling asymmetric vortex sheet pairs. Section \ref{sec:exist} is devoted to the existence of  co-rotating asymmetric vortex sheets. We introduce the nonlinear functional, then we prove its continuity and its $C^1$ regularity, we also compute its linearization at the point vortex pair, we determine the angular velocity and the centroid, and we apply the implicit function theorem. Finally, Section \ref{sec4} deals with the existence of  traveling asymmetric vortex sheets.

\section{Preliminaries}\label{sec2}

\subsection{Notation and function spaces}

Throughout the paper we use the mean value
\begin{equation*}
\fint p(\tau)\dtau:=\frac{1}{2\pi}\int_{0}^{2\pi}p(\tau)\dtau,
\end{equation*}
and the projection onto the mean $\Pm p:=\fint p(\tau)\dtau$. We write $\bfe(x):=(\cos x,\sin x)$, $\bfe'(x)=(-\sin x,\cos x)$, $\bfe_1=(1,0)$, $\bfe_2=(0,1)$, and $Q_\theta$ denotes the counterclockwise rotation of angle $\theta$ about the origin. For a general function $u$ we denote
\begin{equation*}
\Delta u:=u(x)-u(y),\qquad u:=u(x),\qquad \tilde u:=u(y).
\end{equation*}
The identities relating with the singular integrals appearing in this paper, are recalled in Appendix \ref{app:hilbert}, see \eqref{HL}--\eqref{ident}.

Let $a>0$ be fixed and let $S_a:=\{z\in\C:\ |\mathrm{Im}(z)|\le a\}$. Let $C_w(a)$  be  the space of $2\pi$-periodic functions which are holomorphic in the interior of $S_a$, continuous in $S_a$ and real on $\R$. For $k\in\N$ we define the even and odd spaces
\begin{align*}
X^k&:=\Big\{u\in C_w(a):\ u(x)=\sum_{j=1}^\infty a_j\cos(jx),\ \ \norm{u}_{X^k}:=\sum_{\pm}\norm{u(\cdot\pm a\mathbf i)}_{H^k(\T)}<\infty\Big\},\\
Y^k&:=\Big\{u\in C_w(a):\ u(x)=\sum_{j=1}^\infty a_j\sin(jx),\ \ \norm{u}_{Y^k}:=\sum_{\pm}\norm{u(\cdot\pm a\mathbf i)}_{H^k(\T)}<\infty\Big\},
\end{align*}
and,  the space with a restriction on the first Fourier coefficient
\begin{equation}\label{Ak}
\mathcal A^k:=\Big\{(r,s)\in X^{k+1}\times X^k:\ \fint r(x)\cos(x)\dx=-\fint s(x)\cos(x)\dx\Big\},
\end{equation}
endowed with the norm $\norm{(r,s)}_{\mathcal A^k}:=\norm{r}_{X^{k+1}}+\norm{s}_{X^k}$. In other words, writing $r=\displaystyle\sum_{j}a_j\cos(jx)$ and $s=\displaystyle\sum_{j}b_j\cos(jx)$, we have that $(r,s)\in\mathcal A^k$ if and only if $b_1=-a_1$.

We also introduce the space
\begin{equation}\label{Bk}
\mathcal B^k:=\Big\{(u,v)\in Y^k\times X^k:\ \fint u(x)\sin x\dx=-\fint v(x)\cos x\dx\Big\},
\end{equation}
endowed with the norm $\norm{u}_{Y^k}+\norm{v}_{X^k}$. In other words, writing $u=\displaystyle\sum_{j}A_j\sin(jx)$ and $v=\displaystyle\sum_jB_j\cos(jx)$, we have that $(u,v)\in\mathcal B^k$ if and only if $B_1=-A_1$. In contrast with \cite{CQZ}, the space $\mathcal B^k$ is not used in the proofs of our results, since we take $\Omega$ and $x_0$ (respectively $U$ and $\kappa_2$) as unknowns, the functionals take values in the whole space $Y^k\times X^k\times Y^k\times X^k$. 

To eliminate possible singularities at $\ep=0$, we will frequently use the Taylor formula 
\begin{equation}\label{taylor}
\frac{1}{(A+B)^{\lambda}}=\frac{1}{A^{\lambda}}-\lambda\int_0^1\frac{B}{(A+tB)^{1+\lambda}}\dt .
\end{equation}

Finally, the existence of solutions will be obtained from the implicit function theorem, which we recall in Appendix \ref{app:ift}.

\subsection{ Co-rotating asymmetric vortex sheets}\label{sec:rot}

Let $\ep>0$ be small, $d>0$, $b_1,b_2>0$ and $\kappa_1,\kappa_2\in\R\setminus\{0\}$. We consider an initial vorticity of the form
\begin{equation}\label{sole}
\vartheta_{0,\ep}=\frac{\kappa_1}{\ep b_1}\omega_1(t,x)\delta_{\Gamma^\ep_1}+\frac{\kappa_2}{\ep b_2}\omega_2(t,x)\delta_{\Gamma^\ep_2},
\end{equation}
where $\Gamma_i^\ep=\{z_i(t,x):\ x\in[0,2\pi)\}$, $i=1,2$, are two simple closed curves, $\omega_i(t,x)$ is the vortex strength with respect to the parametrization $z_i$, and $\delta_{\Gamma}$ is the Dirac measure on the curve $\Gamma$. The curves $\Gamma_1^\ep$ and $\Gamma_2^\ep$ will be small perturbations of the circles of radii $\ep b_1$ and $\ep b_2$ centered at $0$ and $d\bfe_1$. 
\begin{equation}\label{BRsol}
\vv(t,z_i(t,x))=-\frac{1}{2\pi}\sum_{j=1}^2\kappa_j\,PV\!\!\int_0^{2\pi}\frac{(z_i(t,x)-z_j(t,y))^\perp}{|z_i(t,x)-z_j(t,y)|^2}\,\gamma_j(t,y)\dy,\qquad i=1,2,
\end{equation}
while, since $\frac{\kappa_i}{\ep b_i}$ is a constant, \eqref{1-4} is equivalent to the same equation for $\omega_i$.

We look for solutions rotating uniformly about the point $(x_0,0)$ with angular velocity $\Omega$, namely
\begin{equation}\label{rotsol}
\vartheta_\ep(t,X)=\vartheta_{0,\ep}\big(Q_{-\Omega t}(X-x_0\bfe_1)+x_0\bfe_1\big),
\end{equation}
where $\Gamma_i^\ep(t)-x_0\bfe_1=Q_{\Omega t}\big(\Gamma_i^\ep-x_0\bfe_1\big)$. In the frame rotating with angular velocity $\Omega$ about $(x_0,0)$, the evolution equations \eqref{1-3} and \eqref{1-4} become
\begin{align}
\partial_tz_i(t,x)&=\vv(t,z_i)+\Omega\big(z_i(t,x)-x_0\bfe_1\big)^\perp+c_i(t,x)\,\partial_xz_i(t,x),\label{rotframe}\\
\partial_t\omega_i(t,x)&=\partial_x\big(c_i(t,x)\omega_i(t,x)\big),\label{rotframe2}
\end{align}
with $c_i(t,x)$ being the reparametrization freedom. Since $z_i(t,x)$ parametrizes the same curve as $z_i(0,x)$, $\partial_tz_i(t,x)$ is tangent to $\Gamma_i^\ep$. Let $\mathbf n(z_i(t,x))$ be the unit outward normal vector at $z_i(t,x)$, and $\mathbf s(z_i(t,x))$ be the unit tangent vector oriented counterclockwise. Multiplying \eqref{rotframe} by $\mathbf n(z_i(t,x))$ and using $\mathbf n(z_i(t,x))\cdot\partial_xz_i(t,x)=0$, we obtain
\begin{equation}\label{normal}
\Big(\vv(t,z_i)+\Omega\big(z_i(t,x)-x_0\bfe_1\big)^\perp\Big)\cdot\mathbf n(z_i(t,x))=0,\qquad\forall\,x\in[0,2\pi).
\end{equation}
To eliminate the time parameter and derive another necessary relation, we choose
\begin{equation*}
c_i(t,x)=-\frac{\Big(\vv(t,z_i)+\Omega\big(z_i(t,x)-x_0\bfe_1\big)^\perp\Big)\cdot\mathbf s(z_i(t,x))}{|\partial_xz_i(t,x)|},
\end{equation*}
and hence obtain $\partial_tz_i(t,x)\cdot\mathbf s(z_i(t,x))=0$ by multiplying \eqref{rotframe} by $\mathbf s(z_i(t,x))$. Since $\mathbf s(z_i(t,x))$ is parallel to $\partial_tz_i(t,x)$ and has unit length, it holds $\partial_tz_i(t,x)=0$, which implies that $z_i(t,x)=z_i(x)$ remains fixed in time. On the other hand, $\vartheta_\ep$ is stationary in the rotating frame. Since $z_i$ is invariant, the vortex strength will not evolve with time, namely $\omega_i(t,x)=\omega_i(x)$, and hence also $\gamma_i(t,x)=\gamma_i(x)$. From \eqref{rotframe2} we derive $c_i(x)\omega_i(x)=C_i$ for some constant $C_i$, and, according to the definition of $c_i$, we have
\begin{equation}\label{tangential}
\frac{\Big(\vv(z_i)+\Omega\big(z_i(x)-x_0\bfe_1\big)^\perp\Big)\cdot\mathbf s(z_i(x))\,\omega_i(x)}{|\partial_xz_i(x)|}=-C_i,\qquad\forall\,x\in[0,2\pi).
\end{equation}
Conversely, if \eqref{normal} and \eqref{tangential} hold for $i=1,2$, then $\vartheta_\ep$ defined in \eqref{rotsol} is a solution of \eqref{1-1}, with the above choice of $c_i$.

Now we can assume that the curves are parametrized as
\begin{equation}\label{curves}
\begin{gathered}
z_1(x)=\big(\ep b_1R_1(x)\cos(x),\ \ep b_1R_1(x)\sin(x)\big),\\
z_2(x)=\big(d-\ep b_2R_2(x)\cos(x),\ -\ep b_2R_2(x)\sin(x)\big),
\end{gathered}
\qquad x\in[0,2\pi),
\end{equation}
for some functions $R_1(x)$, $R_2(x)$ close to $1$. The second curve is the image of a curve centered at the origin under the rotation of angle $\pi$ about $(\frac d2,0)$. This makes the system for $i=2$ formally identical to the system for $i=1$. Both parametrizations are counterclockwise, so that $\mathbf n(z_i)=\frac{(\partial_xz_i)^\perp}{|\partial_xz_i|}$ and $\mathbf s(z_i)=\frac{\partial_xz_i}{|\partial_xz_i|}$. Writing $\bfe(x)=(\cos x,\sin x)$, we have $z_1(x)=\ep b_1R_1(x)\bfe(x)$, $z_2(x)=d\bfe_1-\ep b_2R_2(x)\bfe(x)$. 
We set
\begin{equation}\label{xbi}
x_{0,1}:=x_0,\qquad x_{0,2}:=d-x_0,\qquad\text{i.e.}\quad x_{0,i}=(-1)^{i+1}x_0+(i-1)d,
\end{equation}
which is the distance from the center of the $i$-th circle to the center of rotation. We see that constructing co-rotating asymmetric vortex sheet pairs is equivalent to finding a zero of the following equations, for $i=1,2$,
\begin{equation}\label{rot}
\begin{split}
&\Omega\Big(\ep b_iR_i(x)R_i'(x)-x_{0,i}\big(R_i'(x)\cos x-R_i(x)\sin x\big)\Big)\\
&{\resizebox{.95\hsize}{!}{$-\dfrac{\kappa_i}{\ep b_i}\ds\fint\frac{R_i(x)R_i'(x)-R_i'(x)R_i(y)\cos(x-y)+R_i(x)R_i(y)\sin(x-y)}{\left(R_i(x)-R_i(y)\right)^2+4R_i(x)R_i(y)\sin^2\left(\frac{x-y}{2}\right)}\,\gamma_i(y)\dy$}}\\
&{\resizebox{.95\hsize}{!}{$-\kappa_{3-i}\ds\fint\frac{\ep b_iR_i(x)R_i'(x)-d\big(R_i'(x)\cos x-R_i(x)\sin x\big)+\ep b_{3-i}R_{3-i}(y)\big(R_i'(x)\cos(x-y)-R_i(x)\sin(x-y)\big)}{\Dinter{i}}\,\gamma_{3-i}(y)\dy$}}\\
&\qquad=0,
\end{split}
\end{equation}
and to construct traveling asymmetric vortex sheet pairs is equivalent to finding a zero of the following equations, for $i=1,2$,
\begin{equation}\label{rotT}
\begin{split}
&(I-\Pm)\Bigg\{\frac{\gamma_i(x)}{R_i'(x)^2+R_i(x)^2}\Bigg[\Omega\Big(-\ep b_iR_i(x)^2+x_{0,i}\big(R_i'(x)\sin x+R_i(x)\cos x\big)\Big)\\
&{\resizebox{.95\hsize}{!}{$-\dfrac{\kappa_i}{\ep b_i}\ds\fint\frac{-R_i(x)^2+R_i(x)R_i(y)\cos(x-y)+R_i'(x)R_i(y)\sin(x-y)}{\left(R_i(x)-R_i(y)\right)^2+4R_i(x)R_i(y)\sin^2\left(\frac{x-y}{2}\right)}\,\gamma_i(y)\dy$}}\\
&{\resizebox{.95\hsize}{!}{$+\kappa_{3-i}\ds\fint\frac{\ep b_iR_i(x)^2-d\big(R_i'(x)\sin x+R_i(x)\cos x\big)+\ep b_{3-i}R_{3-i}(y)\big(R_i'(x)\sin(x-y)+R_i(x)\cos(x-y)\big)}{\Dinter{i}}\,\gamma_{3-i}(y)\dy\Bigg]\Bigg\}$}}\\
&\qquad=0.
\end{split}
\end{equation}
The constant $-C_i$ is removed by the projection $I-\Pm$. 
\subsection{ Traveling asymmetric vortex sheets}\label{sec:trav}

For the traveling case,  we consider vortex sheets of opposite signs,
\begin{equation}\label{initial}
\vartheta_{0,\ep}=\frac{\kappa_1}{\ep b_1}\omega_1(t,x)\delta_{\Gamma^\ep_1}-\frac{\kappa_2}{\ep b_2}\omega_2(t,x)\delta_{\Gamma^\ep_2}.
\end{equation}
We look for solutions traveling with constant speed $U$ in the direction $\bfe_2$,
\begin{equation}\label{travsol}
\vartheta_\ep(t,X)=\vartheta_{0,\ep}(X-tU\bfe_2).
\end{equation}
By a similar method as in Section \ref{sec:rot}, we can make $\vartheta_\ep$ stationary in the translating frame, and $z_i$, $\gamma_i$ independent of time. Hence the two equations for traveling vortex sheet pairs are, for $i=1,2$,
\begin{equation}\label{travNT}
\big(\vv(z_i)-U\bfe_2\big)\cdot\mathbf n(z_i(x))=0,\qquad \frac{\big(\vv(z_i)-U\bfe_2\big)\cdot\mathbf s(z_i(x))\,\omega_i(x)}{|\partial_xz_i(x)|}=-C_i,\qquad\forall\,x\in[0,2\pi),
\end{equation}
We see that constructing traveling asymmetric vortex sheet pairs is equivalent to finding a zero of the following equations, for $i=1,2$,
\begin{equation}\label{t3}
\begin{split}
&-U\big(R_i(x)\sin x-R_i'(x)\cos x\big)\\
&-\dfrac{\kappa_i}{\ep b_i}\ds\fint\frac{R_i(x)R_i'(x)-R_i'(x)R_i(y)\cos(x-y)+R_i(x)R_i(y)\sin(x-y)}{\left(R_i(x)-R_i(y)\right)^2+4R_i(x)R_i(y)\sin^2\left(\frac{x-y}{2}\right)}\,\gamma_i(y)\dy\\
&{\resizebox{.95\hsize}{!}{$+\kappa_{3-i}\ds\fint\frac{\ep b_iR_i(x)R_i'(x)-d\big(R_i'(x)\cos x-R_i(x)\sin x\big)+\ep b_{3-i}R_{3-i}(y)\big(R_i'(x)\cos(x-y)-R_i(x)\sin(x-y)\big)}{\Dinter{i}}\,\gamma_{3-i}(y)\dy$}}\\
&\qquad=0,
\end{split}
\end{equation}
and
\begin{equation}\label{t3T}
\begin{split}
&(I-\Pm)\Bigg\{\frac{\gamma_i(x)}{R_i'(x)^2+R_i(x)^2}\Bigg[-U\big(R_i'(x)\sin x+R_i(x)\cos x\big)\\
&-\dfrac{\kappa_i}{\ep b_i}\ds\fint\frac{-R_i(x)^2+R_i(x)R_i(y)\cos(x-y)+R_i'(x)R_i(y)\sin(x-y)}{\left(R_i(x)-R_i(y)\right)^2+4R_i(x)R_i(y)\sin^2\left(\frac{x-y}{2}\right)}\,\gamma_i(y)\dy\\
&{\resizebox{.95\hsize}{!}{$-\kappa_{3-i}\ds\fint\frac{\ep b_iR_i(x)^2+\ep b_{3-i}R_{3-i}(y)\big(R_i'(x)\sin(x-y)+R_i(x)\cos(x-y)\big)}{\Dinter{i}}\,\gamma_{3-i}(y)\dy$}}\\
&
{\resizebox{.95\hsize}{!}{$+\kappa_{3-i}\ds\fint\frac{d\big(R_i'(x)\sin x+R_i(x)\cos x\big)}{\Dinter{i}}\,\gamma_{3-i}(y)\dy\Bigg]\Bigg\}$}}\\
&\qquad=0.
\end{split}
\end{equation}
Compared with \eqref{rot}--\eqref{rotT}, the rotation terms are replaced by the translation terms, and the interaction terms change sign because the two sheets have circulations of opposite signs.

\section{Existence of asymmetric co-rotating vortex sheet pairs}\label{sec:exist}

In this section, we extend the range of $\ep$ to a neighborhood of $0$ in the case of asymmetric vortex sheet pairs which are small perturbations of a pair of point vortices. In this case, we can assume that the functions $R_i$ and the vortex strengths $\gamma_i$  are of the form
\begin{equation}\label{ansatz}
R_i(x)=1+\ep b_ir_i(x),\qquad \gamma_i(x)=1+\ep b_is_i(x),\qquad (r_i,s_i)\in\mathcal A^k,\quad i=1,2 .
\end{equation}
Now, we describe the set of solutions of the system $\mathcal F=0$ around the point $(\ep,\Omega,x_0,r_1,s_1,r_2,s_2)$. To do this we consider the functionals
\begin{equation*}
\fit[\hsize]{\mathcal F_i=\mathcal F_{i1}+\mathcal F_{i2}+\mathcal F_{i3},\qquad \mathcal H_i=\mathcal H_{i1}+\mathcal H_{i2}+\mathcal H_{i3},\qquad \mathcal H_{il}:=(I-\Pm)\tilde{\mathcal H}_{il},\quad l=1,2,3,\quad\mbox{for}\ i=1,2,}
\end{equation*}
where the terms are
\begin{equation}\label{F1}
\fit[\hsize]{\mathcal F_{i1}=\Omega\Big(\ep^2b_i^2\big(1+\ep b_ir_i(x)\big)r_i'(x)-\ep b_ix_{0,i}r_i'(x)\cos x+x_{0,i}\big(1+\ep b_ir_i(x)\big)\sin x\Big),}
\end{equation}

\begin{equation}\label{F2}
\begin{split}
\mathcal F_{i2}=&\fit[.85\hsize]{{}-\kappa_i\fint\frac{\big(1+\ep b_ir_i(x)\big)r_i'(x)\big(1-\cos(x-y)\big)\big(1+\ep b_is_i(y)\big)}{\Dself{i}}\dy}\\
&\fit[.85\hsize]{{}-\kappa_i\ep b_i\fint\frac{r_i'(x)\big(r_i(x)-r_i(y)\big)\cos(x-y)\big(1+\ep b_is_i(y)\big)}{\Dself{i}}\dy}\\
&\fit[.85\hsize]{{}+\kappa_i\fint\frac{\big(1+\ep b_ir_i(x)\big)\big(r_i(x)-r_i(y)\big)\sin(x-y)\big(1+\ep b_is_i(y)\big)}{\Dself{i}}\dy}\\
&\fit[.85\hsize]{{}-\frac{\kappa_i}{\ep b_i}\fint\frac{\big(1+\ep b_ir_i(x)\big)^2\sin(x-y)\big(1+\ep b_is_i(y)\big)}{\Dself{i}}\dy}\\
=&\ \mathcal F_{i21}+\mathcal F_{i22}+\mathcal F_{i23}+\mathcal F_{i24},
\end{split}
\end{equation}

\begin{equation}\label{F3}
\begin{split}
\mathcal F_{i3}=
&\kappa_{3-i}\ep b_{3-i}\\
&\fit[.85\hsize]{\qquad\times\fint\frac{\big(1+\ep b_ir_i(x)\big)\big(1+\ep b_{3-i}r_{3-i}(y)\big)\sin(x-y)\big(1+\ep b_{3-i}s_{3-i}(y)\big)}{\Dinter{i}}\dy}\\
&-\kappa_{3-i}\ep b_i\\
&\fit[.85\hsize]{\qquad\times\fint\frac{r_i'(x)\big(\ep b_i(1+\ep b_ir_i(x))+\ep b_{3-i}(1+\ep b_{3-i}r_{3-i}(y))\cos(x-y)\big)\big(1+\ep b_{3-i}s_{3-i}(y)\big)}{\Dinter{i}}\dy},
\end{split}
\end{equation}

\begin{equation}\label{F4}
\begin{split}
\mathcal F_{i4}=  &  \kappa_{3-i}\fit[.85\hsize]{\fint\frac{\,\ep b_ir_i'(x)d\cos x(1+\ep b_{3-i}s_{3-i}(y))}{\Dinter{i}}\dy}\\
&
 -\kappa_{3-i}\fit[.85\hsize]{\fint\frac{\,d\big(1+\ep b_ir_i(x)\big)\sin x(1+\ep b_{3-i}s_{3-i}(y))}{\Dinter{i}}\dy}
\end{split}
\end{equation}
and

\begin{equation}\label{H1}
\fit[\hsize]{\tilde{\mathcal H}_{i1}=\Omega\,\frac{1+\ep b_is_i(x)}{\Ufull{i}}\Big(-\ep b_i\big(1+\ep b_ir_i(x)\big)^2+x_{0,i}\big(\ep b_ir_i'(x)\sin x+(1+\ep b_ir_i(x))\cos x\big)\Big),}
\end{equation}

\begin{equation}\label{H2}
\begin{split}
\tilde{\mathcal H}_{i2}=&\fit[.85\hsize]{{}-\frac{\kappa_i}{\ep b_i}\,\frac{1+\ep b_is_i(x)}{\Ufull{i}}}\\
&\fit[.85\hsize]{\qquad\times\fint\frac{\big(1+\ep b_ir_i(x)\big)\big(1+\ep b_ir_i(y)\big)\big(\cos(x-y)-1\big)\big(1+\ep b_is_i(y)\big)}{\Dself{i}}\dy}\\
&\fit[.85\hsize]{{}+\kappa_i\ep b_i\,\frac{1+\ep b_is_i(x)}{\Ufull{i}}}\\
&\fit[.85\hsize]{\qquad\times\fint\frac{r_i'(x)\big(r_i(x)-r_i(y)\big)\sin(x-y)\big(1+\ep b_is_i(y)\big)}{\Dself{i}}\dy}\\
&\fit[.85\hsize]{{}+\kappa_i\,\frac{1+\ep b_is_i(x)}{\Ufull{i}}}\\
&\fit[.85\hsize]{\qquad\times\fint\frac{\big(1+\ep b_ir_i(x)\big)\big(r_i(x)-r_i(y)\big)\big(1+\ep b_is_i(y)\big)}{\Dself{i}}\dy}\\
&\fit[.85\hsize]{{}-\kappa_i\,\frac{1+\ep b_is_i(x)}{\Ufull{i}}}\\
&\fit[.85\hsize]{\qquad\times\fint\frac{r_i'(x)\big(1+\ep b_ir_i(x)\big)\sin(x-y)\big(1+\ep b_is_i(y)\big)}{\Dself{i}}\dy}\\
=&\ \tilde{\mathcal H}_{i21}+\tilde{\mathcal H}_{i22}+\tilde{\mathcal H}_{i23}+\tilde{\mathcal H}_{i24},
\end{split}
\end{equation}

\begin{equation}\label{H3}
\begin{split}
\tilde{\mathcal H}_{i3}=&\fit[.85\hsize]{{}-\kappa_{3-i}\,d\,\frac{1+\ep b_is_i(x)}{\Ufull{i}}\big(1+\ep b_ir_i(x)\big)\cos x}\\
&\fit[.85\hsize]{\qquad\times\fint\frac{1+\ep b_{3-i}s_{3-i}(y)}{\Dinter{i}}\dy}\\
&\fit[.85\hsize]{{}+\kappa_{3-i}\frac{1+\ep b_is_i(x)}{\Ufull{i}}\big(1+\ep b_ir_i(x)\big)}\\
&\fit[.85\hsize]{\qquad\times\fint\frac{\big(\ep b_i(1+\ep b_ir_i(x))+\ep b_{3-i}(1+\ep b_{3-i}r_{3-i}(y))\cos(x-y)\big)\big(1+\ep b_{3-i}s_{3-i}(y)\big)}{\Dinter{i}}\dy}\\
=&\ \tilde{\mathcal H}_{i31}+\tilde{\mathcal H}_{i32}.
\end{split}
\end{equation}
\begin{equation}
     \begin{split}
  \tilde{\mathcal H}_{i4}     =  &\fit[.85\hsize]{{}\kappa_{3-i}\ep b_ir_i'(x)\frac{1+\ep b_is_i(x)}{\Ufull{i}}}\\
&\fit[.85\hsize]{\qquad\times\fint\frac{\ep b_{3-i}(1+\ep b_{3-i}r_{3-i}(y))\sin(x-y)\big(1+\ep b_{3-i}s_{3-i}(y)\big)}{\Dinter{i}}\dy}\\
&-\fit[.85\hsize]{{}\kappa_{3-i}\ep b_ir_i'(x)\frac{1+\ep b_is_i(x)}{\Ufull{i}}}\\
&\fit[.85\hsize]{\qquad\times\fint\frac{d\sin x\big(1+\ep b_{3-i}s_{3-i}(y)\big)}{\Dinter{i}}\dy}\\
=&\ \tilde{\mathcal H}_{i41}+\tilde{\mathcal H}_{i42}.
\end{split}
\end{equation}
We also denote $\mathcal H_{i2l}:=(I-\Pm)\tilde{\mathcal H}_{i2l}$, $l=1,2,3,4$, and $\mathcal H_{i3l}:=(I-\Pm)\tilde{\mathcal H}_{i3l}$, $l=1,2,3$.

%%%
%%%
%%%

\vspace{3mm}

Now, we outline the proof of our main results. To do this, we shall use the following notation related to the nonlinear functionals
\begin{equation*}
\begin{gathered}
\bfF(\ep,\Omega,x_0,r_1,s_1,r_2,s_2):=\big(\mathcal F_1,\mathcal H_1,\mathcal F_2,\mathcal H_2\big)(\ep,\Omega,x_0,r_1,s_1,r_2,s_2),\\
\bfG(\ep,U,\kappa_2,r_1,s_1,r_2,s_2):=\big(\mathcal G_1,\mathcal J_1,\mathcal G_2,\mathcal J_2\big)(\ep,U,\kappa_2,r_1,s_1,r_2,s_2),
\end{gathered}
\end{equation*}
where $\mathcal F_i,\mathcal G_i$ and $\mathcal H_i,\mathcal J_i$ are the normal and tangential equations of the $i$-th sheet, see Sections \ref{sec:exist} and \ref{sec4}. From now on, we define $\gfr=(\Omega,x_0,r_1,s_1,r_2,s_2)$. In order to ensure the existence of co-rotating asymmetric vortex sheet pairs we split the proof into several steps. These steps are given in the form of lemmas, and they are dedicated to verifying that the hypotheses of the implicit function theorem for the functional $\bfF(\ep,\gfr)=0$ are fulfilled. More precisely, we need to verify that $\bfF(\ep,\gfr)$ and its Gateaux derivatives in the directions $(h_i,\eta_i)$ and $(h_{3-i},\eta_{3-i})$, $i=1,2$, are continuous from $(-\ep_1,\ep_1)\times\R\times\R\times B_\varrho$ to $\mathcal Y^k:=Y^k\times X^k\times Y^k\times X^k$, where $B_\varrho$ is a neighborhood of the origin in $\mathcal A^k\times\mathcal A^k$, see \eqref{Ak}. Here the main point is that the terms of order $\frac{1}{\ep}$ in the self-induced terms cancel, see Lemmas \ref{lem:cont} and \ref{lem:C1}.

To this end, we need to compute the linearized operator, which is given by
\begin{equation*}
D_{(r_1,s_1,r_2,s_2)}\bfF(0,\Omega,x_0,0,0,0,0)(h_1,\eta_1,h_2,\eta_2)(x)=\sum_{j=1}^\infty\begin{pmatrix}
\frac{\kappa_1}{2}\big(ja_j^1-b_j^1\big)\sin(jx)\\[2pt]
\frac{\kappa_1}{2}\big((j-2)a_j^1+b_j^1\big)\cos(jx)\\[2pt]
\frac{\kappa_2}{2}\big(ja_j^2-b_j^2\big)\sin(jx)\\[2pt]
\frac{\kappa_2}{2}\big((j-2)a_j^2+b_j^2\big)\cos(jx)
\end{pmatrix},
\end{equation*}
where $h_i=\displaystyle\sum_{j\ge1}a_j^i\cos(jx)$ and $\eta_i=\displaystyle\sum_{j\ge1}b_j^i\cos(jx)$. This operator is not invertible from $\mathcal A^k\times\mathcal A^k$ to $\mathcal Y^k$, since for $j=1$ the matrix $\frac{\kappa_i}{2}\begin{pmatrix}1&-1\\-1&1\end{pmatrix}$ is singular. Since $b_1^i=-a_1^i$ in $\mathcal A^k$, the first Fourier mode of the image is $\kappa_ia_1^i(\sin x,-\cos x)$ for the $i$-th sheet. Hence the operator is invertible from $\mathcal A^k\times\mathcal A^k$ onto $\mathcal B^k\times\mathcal B^k$, where $\mathcal B^k$ is defined in \eqref{Bk}. But the direction $(\sin x,\cos x)$ of each sheet is not in its range. One way to deal with this defect is to consider the operator $D_\gfr\bfF(0,\gfr_0)\hfr(x)$, where $\hfr=(\beta_1,\beta_2,h_1,\eta_1,h_2,\eta_2)\in\mathcal X^k:=\R\times\R\times\mathcal A^k\times\mathcal A^k$ and $\gfr_0=(\Omega^*,x_0^*,0,0,0,0)$, satisfying
\begin{equation*}
\begin{split}
D_\gfr\bfF(0,\gfr_0)\hfr(x)={}&\beta_1\frac{d}{\kappa_1+\kappa_2}\begin{pmatrix}\kappa_2\sin x\\ \kappa_2\cos x\\ \kappa_1\sin x\\ \kappa_1\cos x\end{pmatrix}+\beta_2\frac{\kappa_1+\kappa_2}{d^2}\begin{pmatrix}\sin x\\ \cos x\\ -\sin x\\ -\cos x\end{pmatrix}\\
&\qquad+\sum_{j=1}^\infty\begin{pmatrix}
\frac{\kappa_1}{2}\big(ja_j^1-b_j^1\big)\sin(jx)\\[2pt]
\frac{\kappa_1}{2}\big((j-2)a_j^1+b_j^1\big)\cos(jx)\\[2pt]
\frac{\kappa_2}{2}\big(ja_j^2-b_j^2\big)\sin(jx)\\[2pt]
\frac{\kappa_2}{2}\big((j-2)a_j^2+b_j^2\big)\cos(jx)
\end{pmatrix}.
\end{split}
\end{equation*}
The first two terms provide the missing directions $(\sin x,\cos x)$, one for each sheet. They span both directions if and only if
\begin{equation*}
\det\begin{pmatrix}\frac{d\kappa_2}{\kappa_1+\kappa_2} & \frac{\kappa_1+\kappa_2}{d^2}\\[4pt] \frac{d\kappa_1}{\kappa_1+\kappa_2} & -\frac{\kappa_1+\kappa_2}{d^2}\end{pmatrix}=-\frac{\kappa_1+\kappa_2}{d}\neq0 .
\end{equation*}
Therefore, if $\kappa_1+\kappa_2\neq0$, one verifies that $D_\gfr\bfF(0,\gfr_0)$ is an isomorphism from $\mathcal X^k$ to $\mathcal Y^k$, see Lemma \ref{lem:iso}. Hence, we are in a position to invoke the implicit function theorem around the point $(\ep,\Omega,x_0,r_1,s_1,r_2,s_2)=(0,\Omega^*,x_0^*,0,0,0,0)$. This guarantees the existence of non-trivial co-rotating asymmetric vortex sheet pairs, where $\Omega^*$ and $x_0^*$ are defined in \eqref{ang}.

Similarly, in the traveling case we set $\gfr=(U,\kappa_2,r_1,s_1,r_2,s_2)$ and $\gfr_0=(U^*,\kappa_1,0,0,0,0)$, with $U^*$ given in \eqref{tra}. Then
\begin{equation*}
\begin{split}
D_\gfr\bfG(0,\gfr_0)\hfr(x)={}&-\beta_1\begin{pmatrix}\sin x\\ \cos x\\ \sin x\\ \cos x\end{pmatrix}+\frac{\beta_2}{d}\begin{pmatrix}\sin x\\ \cos x\\ 0\\ 0\end{pmatrix}+\sum_{j=1}^\infty\begin{pmatrix}
\frac{\kappa_1}{2}\big(ja_j^1-b_j^1\big)\sin(jx)\\[2pt]
\frac{\kappa_1}{2}\big((j-2)a_j^1+b_j^1\big)\cos(jx)\\[2pt]
\frac{\kappa_1}{2}\big(ja_j^2-b_j^2\big)\sin(jx)\\[2pt]
\frac{\kappa_1}{2}\big((j-2)a_j^2+b_j^2\big)\cos(jx)
\end{pmatrix},
\end{split}
\end{equation*}
which is an isomorphism from $\mathcal X^k$ to $\mathcal Y^k$ for every $\kappa_1\neq0$, see Lemma \ref{lem:trav}. The implicit function theorem then gives non-trivial traveling asymmetric vortex sheet pairs, where the circulation $\kappa_2=\kappa_2(\ep)$ of the second sheet is determined by the solution.

\subsection{Regularity of the functional}

 For $0<\varrho<1$ and $k\ge3$, let $B_\varrho$ be the open neighborhood of zero
\begin{equation*}
\fit[\hsize]{B_\varrho:=\Big\{(r_1,s_1,r_2,s_2)\in\mathcal A^k\times\mathcal A^k:\ \norm{r_i}_{X^{k+1}}+\norm{s_i}_{X^k}<\varrho,\ i=1,2\Big\}.}
\end{equation*}

In what follows, we obtain the $C^1$-regularity of the nonlinear functional $\mathcal F(\ep,\Omega,x_0,r_1,s_1,r_2,s_2)$. First, we begin by proving the continuity of the functional, and we also see how to extend it to the case $\ep\le0$.

\begin{lemma}\label{lem:cont}
Let $0<\varrho<1$. Then the nonlinear functional
\begin{equation*}
\mathcal F(\ep,\Omega,x_0,r_1,s_1,r_2,s_2):(-\ep_1,\ep_1)\times\R\times\R\times B_\varrho\to\mathcal Y^k
\end{equation*}
is continuous. Moreover, for $i=1,2$,
\begin{gather}
\fit[\hsize]{\mathcal F_{i1}=\Omega x_{0,i}\sin x+\ep\mathcal R_{i1}}\notag,\\
\fit[\hsize]{\mathcal F_{i2}=-\frac{\kappa_i}{2}r_i'+\kappa_i\fint\frac{s_i(y)\sin(x-y)}{4\sin^2(\frac{x-y}{2})}\dy+\ep\mathcal R_{i2},}\label{Fexp}\\
\fit[\hsize]{\mathcal F_{i3}=\ep\mathcal R_{i3}}\notag\\
\fit[\hsize]{\mathcal F_{i4}=-\frac{\kappa_{3-i}}{d}\sin x+\ep\mathcal R_{i4},}\notag\\
\fit[\hsize]{\mathcal H_{i1}=\Omega x_{0,i}\cos x+\ep\mathcal Q_{i1},}\notag\\
\fit[\hsize]{\mathcal H_{i2}=-\kappa_ir_i+\kappa_i\fint\frac{r_i(x)-r_i(y)}{4\sin^2(\frac{x-y}{2})}\dy+\kappa_i\frac12s_i+\ep\mathcal Q_{i2},}\label{Hexp}\\
\mathcal H_{i3}=\ep\mathcal Q_{i3},\notag\\
\mathcal H_{i4}=-\frac{\kappa_{3-i}}{d}\cos x+\ep\mathcal Q_{i3}, \notag
\end{gather}
where $\mathcal R_{il}(\ep,\Omega,x_0,r_1,s_1,r_2,s_2):(-\ep_1,\ep_1)\times\R\times\R\times B_\varrho\to Y^k$ and $\mathcal Q_{il}(\ep,\Omega,x_0,r_1,s_1,r_2,s_2):(-\ep_1,\ep_1)\times\R\times\R\times B_\varrho\to X^k$, $l=1,2,3$, are continuous.
\end{lemma}

\begin{proof}

We denote by $\lot$ the remaining terms, in which every factor carries at most $k-1$ derivatives. We now study the term $\mathcal F_{i1}$. By \eqref{F1}, $\mathcal F_{i1}$ is a polynomial in $\ep,r_i,r_i',\cos x,\sin x$. Since $r_i\in X^{k+1}$, $r_i'\in Y^k$ and $r_i$ is even, then $\mathcal F_{i1}\in Y^k$, and it depends continuously on $(\ep,\Omega,x_0,r_i)$. Moreover, we can rewrite $\mathcal F_{i1}$ as
\begin{equation*}
\fit[\hsize]{\mathcal F_{i1}=\Omega x_{0,i}\sin x+\ep\mathcal R_{i1},\qquad \mathcal R_{i1}:=\Omega b_i\Big(\ep b_i\big(1+\ep b_ir_i(x)\big)r_i'(x)-x_{0,i}\big(r_i'(x)\cos x-r_i(x)\sin x\big)\Big).}
\end{equation*}

Next, we check the parity of the functionals. Notice that $r_j$ and $s_j$ are even functions, $j=1,2$. Then $r_j'$ is odd. Replacing $y$ by $-y$ in \eqref{F2}--\eqref{F4}, one can ensure that $\mathcal F_{i2}$, $\mathcal F_{i3}$ are odd. For instance, the denominator in \eqref{F2} is invariant under $(x,y)\mapsto(-x,-y)$, $r_i(-x)-r_i(-y)=r_i(x)-r_i(y)$, $r_i'(-x)=-r_i'(x)$ and $\sin(-x+y)=-\sin(x-y)$, so each integrand of $\mathcal F_{i2}$ is odd under $(x,y)\mapsto(-x,-y)$.

Now we prove that the range of $\mathcal F_{i2}$ is contained in $Y^k$. To this end, we study the four terms in \eqref{F2} separately.  Note that a potential singularity at $\ep = 0$ can only arise from the zeroth derivative of $\mathcal F_{i2}$, and more specifically, it is restricted to the term $\mathcal F_{i24}$.
In order to isolate the singular part of the kernels, we invoke the Taylor formula \eqref{taylor} with $\lambda=1$, and
\begin{equation*}
\fit[\hsize]{A:=4\sin^2\Big(\frac{x-y}{2}\Big),\qquad B:=\ep^2b_i^2\big(r_i(x)-r_i(y)\big)^2+4\ep b_i\big(r_i(x)+r_i(y)+\ep b_ir_i(x)r_i(y)\big)\sin^2\Big(\frac{x-y}{2}\Big).}
\end{equation*}
Since $\displaystyle\int_0^1\frac{B}{(A+tB)^2}\dt=\frac1A-\frac{1}{A+B}$, we have
\begin{equation}\label{K}
\fit[\hsize]{\frac{1}{\Dself{i}}=\frac{1}{4\sin^2(\frac{x-y}{2})}\Big(1+\ep b_iK_\ep^i(x,y)\Big),}
\end{equation}
where
\begin{equation}\label{Kdef}
\begin{split}
K_\ep^i(x,y)&\fit[.85\hsize]{:=-\frac{\ep b_i\big(r_i(x)-r_i(y)\big)^2+4\big(r_i(x)+r_i(y)+\ep b_ir_i(x)r_i(y)\big)\sin^2(\frac{x-y}{2})}{\Dself{i}}}\\
&\fit[.85\hsize]{{}=-\frac{\ep b_i\Big(\Dq{i}\Big)^2+r_i(x)+r_i(y)+\ep b_ir_i(x)r_i(y)}{\Dn{i}}.}
\end{split}
\end{equation}
An important fact on this decomposition is that $K_\ep^i(x,y)$ is not singular at $x=y$. Because of the factor $\frac{1}{\ep b_i}$, the possible singularity at $\ep=0$ may occur only in $\mathcal F_{i24}$. To prove that this cannot happen, we use the decomposition \eqref{K}. Since the sine function is odd, and using the Taylor formula, we get
\begin{equation}\label{F24}
\begin{split}
\mathcal F_{i24}=&\fit[.85\hsize]{{}-\frac{\kappa_i}{\ep b_i}\fint\frac{\sin(x-y)}{4\sin^2(\frac{x-y}{2})}\dy}-\kappa_i\fint\frac{s_i(y)\sin(x-y)}{4\sin^2(\frac{x-y}{2})}\dy\\
&
\qquad+\kappa_i\fint\frac{(r_i(x)+r_i(y))\sin(x-y)}{4\sin^2(\frac{x-y}{2})}\dy+\ep\mathcal R_{i24}\\
&=-\kappa_i\fint\frac{s_i(y)\sin(x-y)}{4\sin^2(\frac{x-y}{2})}\dy+\kappa_i\fint\frac{(r_i(x)+r_i(y))\sin(x-y)}{4\sin^2(\frac{x-y}{2})}\dy+\ep\mathcal R_{i24} ,
\end{split}
\end{equation}
where $\mathcal R_{i24}$ is not singular at $\ep = 0$. Hence $\mathcal F_{i24}$ is not singular at $\ep=0$.

We begin with the term $\mathcal F_{i21}$. Since $1-\cos(x-y)=2\sin^2(\frac{x-y}{2})$ and making the change of variables $x-y\to y$, we can compute the  $k-th$ derivatives of $\mathcal F_{i21}$ with respect to the spatial variable $x$, we get
\begin{equation*}
\begin{split}
\partial^k \mathcal F_{i21}&=	-	\frac{\kappa_i}{2}\fint(1+\varepsilon b_ir_i(x))\partial^{k+1}r_i(x)(1+\varepsilon b_is_i(x-y))dy\\
&
\qquad -\frac{\kappa_i\varepsilon b_i}{2}\fint(1+\varepsilon b_ir_i(x))r_i'(x)\partial^ks_i(x-y)dy+\textrm{l.o.t.},
\end{split}
\end{equation*}
where, we have used decomposition \eqref{K}, and $\text{l.o.t.}$ denotes lower-order terms. Since $r_i \in X^{k+1}$ and $s_i \in X^k$ with $k \ge 3$, it follows that
\begin{equation*}
\|\partial^i r_i(x \pm a\mathbf{i})\|_{L^\infty} \le C \|r_i\|_{X^{k+1}} < \infty \quad \text{and} \quad \|\partial^i s_i(x \pm a\mathbf{i})\|_{L^\infty} \le C \|s_i\|_{X^k} < \infty
\end{equation*}
for $i = 0, 1, 2$. Applying Hölder's inequality and the Mean Value Theorem, we obtain
\begin{equation*}
\|\partial^k \mathcal{F}_{121}(x \pm a\mathbf{i})\|_{L^2} \le C \|r_i\|_{X^{k+1}} + C \|s_i\|_{X^k} < \infty.
\end{equation*}

Next, we deal with the most singular term $\mathcal F_{i22}$.  Taking $k$ derivatives as before, we find
\begin{equation*}
\begin{split}
\partial^k\mathcal F_{i22}=&\fit[.85\hsize]{{}-\kappa_i\ep b_i\,\partial^{k+1}r_i(x)\fint\frac{\big(r_i(x)-r_i(y)\big)\cos(x-y)\big(1+\ep b_is_i(y)\big)}{\Dself{i}}\dy}\\
&\fit[.85\hsize]{{}-\kappa_i\ep b_i\,r_i'(x)\fint\frac{\big(\partial^kr_i(x)-\partial^kr_i(y)\big)\cos(x-y)\big(1+\ep b_is_i(y)\big)}{\Dself{i}}\dy}\\
&\fit[.85\hsize]{{}-\kappa_i\ep^2b_i^2\,r_i'(x)\fint\frac{\big(r_i(x)-r_i(y)\big)\cos(x-y)\,\partial^ks_i(y)}{\Dself{i}}\dy}\\
&\fit[.85\hsize]{{}+\kappa_i\ep b_i\,r_i'(x)\fint\frac{\big(r_i(x)-r_i(y)\big)\cos(x-y)\big(1+\ep b_is_i(y)\big)}{\big(\Dself{i}\big)^2}}\\
&\fit[.85\hsize]{\qquad\times\Big(\DhA{i}{\partial^kr_i}}\\
&\fit[.85\hsize]{\qquad\quad+\DhB{i}{\partial^kr_i}\Big)\dy+\lot}\\
=&:\mathcal I_{i1}+\mathcal I_{i2}+\mathcal I_{i3}+\mathcal I_{i4}+\lot
\end{split}
\end{equation*}
We estimate each term separately. Using the decomposition \eqref{K} with $\|K^i_\varepsilon(x,y)\|_{L^\infty}\le C(\varepsilon,\|r_i\|_{x^{k+1}})<\infty$ and the fact $|y|^2$ is the equivalent infinitesimal for $4\sin^2\left(\frac{y}{2}\right)$ as $|y|\to 0$. The integral in $\mathcal I_{i1}$ is bounded by 
    \begin{equation*}
    	\begin{split}
    		\| \mathcal I_{i1}(x&\pm a\mathbf{i})\|_{L^2} \\
            &
            \le C\left\|\fint\frac{\kappa_i\varepsilon b_i \partial^{k+1}r_i(x\pm a\mathbf{i})(r_i(x\pm a\mathbf{i})-r_i(x\pm a\mathbf{i}-y))\cos(y)(1+\varepsilon s_i(x\pm a\mathbf{i}-y))}{4\sin^2\left(\frac{y}{2}\right)}dy\right\|_{L^2}\\
    		&\le C\|r_i\|_{X^{k+1}}(1+\|s_i (x\pm a\mathbf{i})\|_{L^\infty})\left|\fint\frac{r_i(x\pm a\mathbf{i})-r_i(x\pm a\mathbf{i}-y)dy}{4\sin^2\left(\frac{y}{2}\right)}\right|\\
            	&= C\|r_i\|_{X^{k+1}}(1+\|s_i (x\pm a\mathbf{i})\|_{L^\infty})\left|\fint\frac{2r_i(x\pm a\mathbf{i})-r_i(x\pm a\mathbf{i}-y)-r_i(x\pm a\mathbf{i}+y)dy}{4\sin^2\left(\frac{y}{2}\right)}\right|\\
    		&\le C\|r_i\|_{X^{k+1}}(1+\|s_i\|_{X^k})\|_{L^\infty})\fint r_i''(x\pm a\mathbf{i})dy\\
    		&\le C\|r_i\|_{X^{k+1}}^2(1+\|s_i\|_{X^k}).
    	\end{split}
    \end{equation*}

By similar estimates as before we can get 
    \begin{equation*}
    	\|\mathcal I_{i2}(x\pm a\mathbf{i})\|_{L^2}\le C\|r_i\|_{X^{k+1}}^2(1+\|s_i\|_{X^k}).
    \end{equation*}
In what follows, we turn our attention to $\mathcal{I}_{i3}$. By Hölder's inequality and the Mean Value Theorem, a direct calculation yields
\begin{equation*}
\begin{split}
    \|\mathcal{I}_{i3}(x \pm a\mathbf{i})\|_{L^2} 
    &\le C \|r_i\|_{X^{k+1}} (1 + \|s_i\|_{X^k}) \left| \fint \frac{r_i(x \pm a\mathbf{i}) - r_i(x \pm a\mathbf{i} - y)}{4\sin^2\left(\frac{y}{2}\right)} dy \right| \\
    &\quad \times \Big( \|r_i'(x \pm a\mathbf{i})\|_{L^\infty} \|r_i''(x \pm a\mathbf{i})\|_{L^\infty} + (1 + \|r_i(x \pm a\mathbf{i})\|_{L^\infty}) \|r_i'(x \pm a\mathbf{i})\|_{L^\infty} \Big) \\
    &\le C \|r_i\|_{X^{k+1}}^3 (1 + \|r_i\|_{X^{k+1}}) (1 + \|s_i\|_{X^k}).
\end{split}
\end{equation*}
Finally, it is straightforward to show that
\begin{equation*}
\|\mathcal{I}_{i4}(x\pm a\mathbf{i})\|_{L^2}\le C\|r_i\|_{X^{k+1}}^2(1+\|s_i\|_{X^k}).
\end{equation*}
Because $\mathcal F_{i22}$ is the most singular term in $\mathcal F_{i2}$, the preceding discussion guarantees that $\|\mathcal F_{i2}\|_{Y^k}\le C\|\mathcal F_{i22}\|_{Y^k}<\infty$. 
From the above estimates and the parity of $\mathcal F_{i2}$, we can conclude that $\norm{\mathcal F_{i2}}_{Y^k}\le C$, that is, the functional $\mathcal F_{i2}$ belongs to $Y^k$.

We now verify the continuity of the functional $\mathcal F$. Once again we only treat the most singular term $\mathcal F_{i22}$ with respecto to $r_i$. The other terms and variables are treated in the same way. Set $\Delta r_i=r_i(x)-r_i(y)$, $r_i=r_i(x)$, $\tilde r_i=r_i(y)$, and
\begin{equation*}
D_\ep(r_i):=\ep^2b_i^2(\Delta r_i)^2+4\big(1+\ep b_ir_i\big)\big(1+\ep b_i\tilde r_i\big)\sin^2\Big(\frac{x-y}{2}\Big).
\end{equation*}
Then for $(r_{i1},s_i),(r_{i2},s_i)$ in $B_\varrho$, we have
\begin{equation*}
\begin{split}
&\mathcal F_{i22}(\ep,r_{i1},s_i)-\mathcal F_{i22}(\ep,r_{i2},s_i)=\\
&\fit[.85\hsize]{\qquad-\kappa_i\ep b_i\big(r_{i1}'-r_{i2}'\big)\fint\frac{\Delta r_{i1}\cos(x-y)\big(1+\ep b_i\tilde s_i\big)}{D_\ep(r_{i1})}\dy}\\
&\fit[.85\hsize]{\qquad-\kappa_i\ep b_ir_{i2}'\fint\frac{(\Delta r_{i1}-\Delta r_{i2})\cos(x-y)\big(1+\ep b_i\tilde s_i\big)}{D_\ep(r_{i1})}\dy}\\
&\fit[.85\hsize]{\qquad-\kappa_i\ep b_ir_{i2}'\fint\Delta r_{i2}\cos(x-y)\big(1+\ep b_i\tilde s_i\big)\Big(\frac{1}{D_\ep(r_{i1})}-\frac{1}{D_\ep(r_{i2})}\Big)\dy}\\
&\qquad=:\ I_1+I_2+I_3 .
\end{split}
\end{equation*}
To estimate $I_3$, we use
\begin{equation}\label{diffD}
\begin{split}
\frac{1}{D_\ep(r_{i1})}-\frac{1}{D_\ep(r_{i2})}&=\frac{D_\ep(r_{i2})-D_\ep(r_{i1})}{D_\ep(r_{i1})D_\ep(r_{i2})},\\
D_\ep(r_{i2})-D_\ep(r_{i1})&=\ep^2b_i^2\big(\Delta r_{i2}-\Delta r_{i1}\big)\big(\Delta r_{i2}+\Delta r_{i1}\big)\\
&\fit[.85\hsize]{\quad+4\ep b_i\Big((r_{i2}-r_{i1})(1+\ep b_i\tilde r_{i2})+(1+\ep b_ir_{i1})(\tilde r_{i2}-\tilde r_{i1})\Big)\sin^2\Big(\frac{x-y}{2}\Big).}
\end{split}
\end{equation}
Since $|D_\ep(r_{i1})D_\ep(r_{i2})|\ge4\sin^4(\frac{x-y}{2})$, each $I_j$ is linear in $r_{i1}-r_{i2}$ and its kernel has the same structure as the kernels of $\mathcal F_{i22}$ studied above. Then,
    \begin{equation*}
    	\begin{split}
    	\|\partial^kI_3(x\pm a\mathbf{i})\|_{L^2}&\le C\|r_{i2}\|_{X^{k+1}}^2\big(\|r_{i1}\|_{X^{k+1}}+\|r_{i2}\|_{X^{k+1}}+C\big)\|r_{i1}-r_{i2}\|_{X^{k+1}}(1+\|s_i\|_{X^k})\\
    	&\le C\|r_{i1}-r_{i2}\|_{X^{k+1}}.
    	\end{split}
    \end{equation*}
Repeating those computations, we obtain
\begin{equation*}
\norm{I_j}_{Y^k}\le C\norm{r_{i1}-r_{i2}}_{X^{k+1}},\qquad j=1,2 .
\end{equation*}

The functional $\mathcal F_{i2}$ is affine in $s_i$, then the continuity with respect to $s_i$ is immediate from the above estimates. Concerning $\ep$, after the decompositions \eqref{F24} and \eqref{H21}, all the kernels are polynomials in $\ep$ divided by powers of $\Dn{i}$, $\Ufull{i}$ and of the denominator of the interaction terms. Thus, the continuity with respect to $\ep$ follows in the same way. With these observations, one can conclude that the nonlinear functional $\mathcal F_{i2}:(-\ep_1,\ep_1)\times\R\times\R\times B_\varrho\to\mathcal Y^k$ is continuous. For further use, we write the functional $\mathcal F_{i2}$ as follows
\begin{equation*}
\begin{split}
\mathcal F_{i21}&=-\kappa_ir_i'(x)\fint\frac{2\sin^2(\frac{x-y}{2})}{4\sin^2(\frac{x-y}{2})}\dy=-\frac{\kappa_i}{2}r_i'+\ep\mathcal R_{i21},\\
\mathcal F_{i22}&=\ep\mathcal R_{i22},\\
\mathcal F_{i23}&=\kappa_i\fint\frac{(r_i(x)-r_i(y))\sin(x-y)}{4\sin^2(\frac{x-y}{2})}\dy+\ep\mathcal R_{i23},\\
\mathcal F_{i24}&=-\kappa_i\fint\frac{s_i(y)\sin(x-y)}{4\sin^2(\frac{x-y}{2})}\dy+\kappa_i\fint\frac{(r_i(x)+r_i(y))\sin(x-y)}{4\sin^2(\frac{x-y}{2})}\dy+\ep\mathcal R_{i24} .
\end{split}
\end{equation*}
Then, one conclude that
\begin{equation}\label{fi2}
\mathcal F_{i2}=-\frac{\kappa_i}{2}\left(r_i'+\displaystyle\fint\frac{s_i(y)\sin(x-y)}{4\sin^2(\frac{x-y}{2})}\dy\right)+\ep\mathcal R_{i2} ,
\end{equation}
where $\mathcal R_{i2}$ is a regular term.

Now, we focus on the functional $\mathcal{H}_{i}$. Similarly as before, we have
$$\mathcal H_{i1}=\Omega x_{0,i}\cos x+\ep\mathcal Q_{i1},$$
where
\begin{equation*}
\begin{split}
&\fit[.85\hsize]{\mathcal Q_{i1}:=\Omega b_i(I-\Pm)\Bigg[-\frac{\big(1+\ep b_is_i(x)\big)\big(1+\ep b_ir_i(x)\big)^2}{\Ufull{i}}}\\
&\qquad \fit[.85\hsize]{\qquad+x_{0,i}\,\frac{\big(1+\ep b_is_i(x)\big)r_i'(x)\sin x+\big(s_i(x)-r_i(x)+\ep b_i(r_i(x)s_i(x)-r_i(x)^2-r_i'(x)^2)\big)\cos x}{\Ufull{i}}\Bigg].}
\end{split}
\end{equation*}
We now prove that the range of $\mathcal H_{i2}$ is contained in $X^k$, proceeding as for $\mathcal F_{i2}$. Notice that the factor $\frac{1}{\ep b_i}$, this is the only term of $\tilde{\mathcal H}_{i2}$ which might be singular at $\ep=0$. Since $\cos(x-y)-1=-2\sin^2(\frac{x-y}{2})$, the decomposition \eqref{K} gives

\begin{equation*}
\begin{split}
\tilde{\mathcal H}_{i21}=&\fit[.85\hsize]{\ \frac{\kappa_i}{2\ep b_i}\,\frac{\big(1+\ep b_is_i(x)\big)\big(1+\ep b_ir_i(x)\big)}{\Ufull{i}}\fint\big(1+\ep b_ir_i(y)\big)\big(1+\ep b_is_i(y)\big)\big(1+\ep b_iK_\ep^i(x,y)\big)\dy}\\
=&\fit[.85\hsize]{\ \frac{\kappa_i}{2\ep b_i}\,\frac{\big(1+\ep b_is_i(x)\big)\big(1+\ep b_ir_i(x)\big)}{\Ufull{i}}}\\
&\fit[.85\hsize]{\times\Big(1+\ep^2b_i^2\fint r_i(y)s_i(y)\dy+\ep b_i\fint\big(1+\ep b_ir_i(y)\big)\big(1+\ep b_is_i(y)\big)K_\ep^i(x,y)\dy\Big),}
\end{split}
\end{equation*}
where we used that $r_i$ and $s_i$ have zero mean. The term of order $\frac{1}{\ep b_i}$ is the constant $\frac{\kappa_i}{2\ep b_i}$, which is eliminated by the projection $I-\Pm$. Therefore
\begin{equation}\label{H21}
\begin{split}
&\fit[.85\hsize]{\mathcal H_{i21}=\frac{\kappa_i}{2}(I-\Pm)\Bigg[\frac{s_i(x)-r_i(x)+\ep b_i\big(r_i(x)s_i(x)-r_i(x)^2-r_i'(x)^2\big)}{\Ufull{i}}}\\
&\fit[.85\hsize]{\qquad\times\Big(1+\ep^2b_i^2\fint r_i(y)s_i(y)\dy+\ep b_i\fint\big(1+\ep b_ir_i(y)\big)\big(1+\ep b_is_i(y)\big)K_\ep^i(x,y)\dy\Big)}\\
&\qquad+\fint\big(1+\ep b_ir_i(y)\big)\big(1+\ep b_is_i(y)\big)K_\ep^i(x,y)\dy\Bigg]
\end{split}
\end{equation}
is not singular at $\ep=0$.  By similar computations as before, we can conclude that 
$$\norm{\mathcal H_{i21}}_{X^k}\le C, \qquad \norm{\mathcal H_{i23}}_{X^k}\le C \qquad \mbox{and}\qquad \norm{\mathcal H_{i24}}_{X^k}\le C .$$

Finally, concerning $\mathcal H_{i24}$, using \eqref{K}, we get
\begin{equation*}
\fit{\tilde{\mathcal H}_{i24}=-\kappa_i\ep b_i\,\frac{(1+\ep b_is_i(x))(1+\ep b_ir_i(x))}{\Ufull{i}}\,r_i'(x)\Big(\fint\frac{s_i(y)\sin(x-y)}{4\sin^2(\frac{x-y}{2})}\dy+\fint\frac{K_\ep^i(x,y)\big(1+\ep b_is_i(y)\big)\sin(x-y)}{4\sin^2(\frac{x-y}{2})}\dy\Big),}
\end{equation*}
and by similar estimates for $\mathcal F_{i24}$, we can conclude that $\norm{\mathcal H_{i24}}_{X^k}\le C$. 

For later use, we write the functional $\mathcal H_{i2}$ as follows, by \eqref{H21}, we have
\begin{equation*}
\fit[\hsize]{\mathcal H_{i21}=\frac{\kappa_i}{2}(I-\Pm)\Big[s_i(x)-r_i(x)-\fint\big(r_i(x)+r_i(y)\big)\dy\Big]=\frac{\kappa_i}{2}\big(s_i-2r_i\big)+\ep\mathcal Q_{i21},}
\end{equation*}
and
\begin{equation*}
\mathcal H_{i22}=\ep\mathcal Q_{i22},\qquad \mathcal H_{i24}=\ep\mathcal Q_{i24},\qquad \mathcal H_{i23}=\kappa_i\fint\frac{r_i(x)-r_i(y)}{4\sin^2(\frac{x-y}{2})}\dy+\ep\mathcal Q_{i23},
\end{equation*}
which gives
\begin{equation}\label{eq:hi2}
\fit[\hsize]{\mathcal H_{i2}=-\kappa_ir_i(x)+\kappa_i\fint\frac{r_i(x)-r_i(y)}{4\sin^2(\frac{x-y}{2})}\dy+\frac{\kappa_i}{2}s_i(x),}
\end{equation}

The functional $\tilde{\mathcal H}_{i2}$ is a polynomial of degree two in $s_i$, so the continuity with respect to $s_i$ is also immediate. Moreover, continuity with respect to $\ep$ follows in the same way. From the above estimates and the parity of $\tilde{\mathcal H}_{i2}$, the functional $\mathcal H_{i2}$ belongs to $X^k$.

We now consider the interaction terms. More precisely, we  study $\mathcal F_{13}$ and $\mathcal F_{14}$. Notice that the denominator of the interaction terms given by $\Dinter{i}$ has positive lower bounds. So $\mathcal F_{i3}$ and $\mathcal F_{i4}$ are less singular than $\mathcal F_{i2}$. Then, $\mathcal F_{i3}$ and $\mathcal F_{i4}$ are both continuous. Using the Taylor's formula \eqref{taylor},

For the interaction terms, to manage the singularity at $\ep=0$, we apply the Taylor formula \eqref{taylor},  we also rewrite $\mathcal F_{i3}$ and $\mathcal F_{i4}$ as
	\begin{equation}\label{2-14}
		\mathcal F_{i3}=\varepsilon \mathcal R_{13},
	\end{equation}
    and
\begin{equation*}
\begin{gathered}
A=d^2,\qquad B=-2d\Big(\ep b_i\big(1+\ep b_ir_i(x)\big)\cos x+\ep b_{3-i}\big(1+\ep b_{3-i}r_{3-i}(y)\big)\cos y\Big)\\
+\Big|\ep b_i\big(1+\ep b_ir_i(x)\big)\bfe(x)+\ep b_{3-i}\big(1+\ep b_{3-i}r_{3-i}(y)\big)\bfe(y)\Big|^2,
\end{gathered}
\end{equation*}
so that $A+B$ is the denominator in \eqref{F3} and $B=O(\ep)$. We arrive at
\begin{equation*}
\begin{split}
\mathcal F_{i41}&\fit[.85\hsize]{{}=-\kappa_{3-i}\,d\big(1+\ep b_ir_i(x)\big)\sin x\fint\big(1+\ep b_{3-i}s_{3-i}(y)\big)\Big(\frac1{d^2}-\int_0^1\frac{B}{(d^2+tB)^2}\dt\Big)\dy}\\
&=-\frac{\kappa_{3-i}}{d}\sin x+\ep\mathcal R_{i41},
\end{split}
\end{equation*}
 where $\mathcal R_{i41}$ is not singular in $\ep$. Hence, 
	\begin{equation}\label{2-15}
		\begin{split}
		\mathcal F_{i4}&=-\frac{\kappa_{3-i}}{d}\sin x+\varepsilon \mathcal R_{i4}
		\end{split}
	\end{equation}
    and in the same way, we have $\mathcal H_{i3}=\ep\mathcal Q_{i3}$ and 
    \begin{equation}\label{eq:hi4}
        \mathcal H_{i4}=-\frac{\kappa_{3-i}}{d}\cos x+\ep\mathcal Q_{i4} ,
    \end{equation}
where $\mathcal R_{i4}$ is not singular in $\ep$.
%%

%%%%%%%%%%%%%%
%%%%%%%%
%%%%%%%%%%%
%%%%%%%%%

%%%
%%

  Together with the expansions of $\mathcal F_{i1}$ and $\mathcal H_{i1}$ obtained above, this proves \eqref{Fexp}--\eqref{Hexp}, and the proof is complete.
\end{proof}

For $(\ep,\Omega,x_0,r_1,s_1,r_2,s_2)\in(-\ep_1,\ep_1)\times\R\times\R\times B_\varrho$ and $(h_i,\eta_i)\in\mathcal A^k$, let
\begin{equation*}
\begin{gathered}
\fit[\hsize]{\partial_{r_i}\mathcal F(\ep,\Omega,x_0,r_1,s_1,r_2,s_2)h_i:=\lim_{t\to0}\frac{1}{t}\Big(\mathcal F(\ep,\Omega,x_0,\ldots,r_i+th_i,\ldots)-\mathcal F(\ep,\Omega,x_0,r_1,s_1,r_2,s_2)\Big),}\\
\fit[\hsize]{\partial_{s_i}\mathcal F(\ep,\Omega,x_0,r_1,s_1,r_2,s_2)\eta_i:=\lim_{t\to0}\frac{1}{t}\Big(\mathcal F(\ep,\Omega,x_0,\ldots,s_i+t\eta_i,\ldots)-\mathcal F(\ep,\Omega,x_0,r_1,s_1,r_2,s_2)\Big)}
\end{gathered}
\end{equation*}
be the Gateaux derivatives of $\mathcal F$ in the directions $h_i$ and $\eta_i$, respectively.  In the following lemma, we compute the Gateaux derivatives of the nonlinear functional $\mathcal F$ and then we prove that they are in fact continuous. In other words, we prove that $\mathcal F$ is a $C^1$ function.

\begin{lemma}\label{lem:C1}
For each $(\ep,\Omega,x_0,r_1,s_1,r_2,s_2)\in(-\ep_1,\ep_1)\times\R\times\R\times B_\varrho$, the partial derivatives of the functional $\mathcal F$ with respect to $\Omega$, $x_0$, $r_i$ and $s_i$, $i=1,2$, exist and are continuous. Moreover, $\mathcal F\in C^1\big((-\ep_1,\ep_1)\times\R\times\R\times B_\varrho;\mathcal Y^k\big)$, and at $\ep=0$ we have, for every $g=(\Omega,x_0,r_1,s_1,r_2,s_2)$ and $i=1,2$,
\begin{equation}\label{lin0}
\begin{gathered}
\fit[\hsize]{\partial_{r_i}\mathcal F_i(0,g)h_i=-\frac{\kappa_i}{2}h_i'(x),\qquad \partial_{s_i}\mathcal F_i(0,g)\eta_i=-\kappa_i\fint\frac{\eta_i(y)\sin(x-y)}{4\sin^2(\frac{x-y}{2})}\dy,}\\
\fit[\hsize]{\partial_{r_i}\mathcal H_i(0,g)h_i=-\kappa_ih_i(x)+\kappa_i\fint\frac{h_i(x)-h_i(y)}{4\sin^2(\frac{x-y}{2})}\dy,\qquad \partial_{s_i}\mathcal H_i(0,g)\eta_i=\frac{\kappa_i}{2}\eta_i(x),}\\
\fit[\hsize]{\partial_{r_{3-i}}\mathcal F_i(0,g)=\partial_{s_{3-i}}\mathcal F_i(0,g)=\partial_{r_{3-i}}\mathcal H_i(0,g)=\partial_{s_{3-i}}\mathcal H_i(0,g)=0,}\\
\partial_\Omega \mathcal F_i(0,g)=x_{0,i}\sin x,\qquad \partial_{x_0}\mathcal F_i(0,g)=(-1)^{i+1}\Omega\sin x,\\
\partial_\Omega \mathcal H_i(0,g)=x_{0,i}\cos x,\qquad \partial_{x_0}\mathcal H_i(0,g)=(-1)^{i+1}\Omega\cos x .
\end{gathered}
\end{equation}
\end{lemma}

\begin{proof}
First of all, notice that the most delicate part is to guarantee the existence and continuity of the Gateaux derivatives in the directions $h_i$ and $\eta_i$. Notice that the functional $\mathcal F$ is linearly dependent on $\Omega$ and affinely dependent on $x_0$ (through $x_{0,i}$ in $\mathcal F_{i1}$ and $\tilde{\mathcal H}_{i1}$), and then it is easy to compute
\begin{equation}\label{dOmega}
\begin{split}
\partial_\Omega \mathcal F_i&\fit[.85\hsize]{{}=\ep^2b_i^2\big(1+\ep b_ir_i(x)\big)r_i'(x)-\ep b_ix_{0,i}r_i'(x)\cos x+x_{0,i}\big(1+\ep b_ir_i(x)\big)\sin x,}\\
\partial_{x_0}\mathcal F_i&=(-1)^{i+1}\Omega\Big(\big(1+\ep b_ir_i(x)\big)\sin x-\ep b_ir_i'(x)\cos x\Big),\\
\partial_\Omega \mathcal H_i&\fit[.85\hsize]{{}=(I-\Pm)\Big[\frac{1+\ep b_is_i(x)}{\Ufull{i}}\Big(-\ep b_i\big(1+\ep b_ir_i(x)\big)^2+x_{0,i}\big(\ep b_ir_i'(x)\sin x+(1+\ep b_ir_i(x))\cos x\big)\Big)\Big],}\\
\partial_{x_0}\mathcal H_i&\fit[.85\hsize]{{}=(-1)^{i+1}\Omega\,(I-\Pm)\Big[\frac{1+\ep b_is_i(x)}{\Ufull{i}}\Big(\ep b_ir_i'(x)\sin x+\big(1+\ep b_ir_i(x)\big)\cos x\Big)\Big],}
\end{split}
\end{equation}
which are clearly continuous, as shown in the study of the rotation terms in the proof of Lemma \ref{lem:cont}.

Now we compute the Gateaux derivatives. From \eqref{F1}, it is easy to see that the first term of the Gateaux derivative with respect to $r_i$ is
\begin{equation}\label{dF1}
\fit[\hsize]{\partial_{r_i}\mathcal F_{i1}h_i=\Omega\ep b_i\Big(\ep b_i\big(1+\ep b_ir_i(x)\big)h_i'(x)+\ep^2b_i^2h_i(x)r_i'(x)-x_{0,i}h_i'(x)\cos x+x_{0,i}h_i(x)\sin x\Big),}
\end{equation}
which is continuous. Similarly, $\partial_{r_{3-i}}\mathcal F_{i1}h_{3-i}=0$ and $\partial_{s_j}\mathcal F_{i1}\eta_j=0$, for $j=1,2$. Next, we claim that $\partial_{r_i}\mathcal F_{i2}h_i=\partial_{r_i}\mathcal F_{i21}h_i+\partial_{r_i}\mathcal F_{i22}h_i+\partial_{r_i}\mathcal F_{i23}h_i+\partial_{r_i}\mathcal F_{i24}h_i$ is continuous, where
\begin{equation}\label{dF21}
\begin{split}
\partial_{r_i}\mathcal F_{i21}h_i=&\fit[.85\hsize]{{}-\kappa_i\fint\frac{\big(\ep b_ih_i(x)r_i'(x)+(1+\ep b_ir_i(x))h_i'(x)\big)\big(1-\cos(x-y)\big)\big(1+\ep b_is_i(y)\big)}{\Dself{i}}\dy}\\
&\fit[.85\hsize]{{}+\kappa_i\fint\frac{\big(1+\ep b_ir_i(x)\big)r_i'(x)\big(1-\cos(x-y)\big)\big(1+\ep b_is_i(y)\big)}{\big(\Dself{i}\big)^2}}\\
&\fit[.85\hsize]{\qquad\times\Big(\DhA{i}{h_i}}\\
&\fit[.85\hsize]{\qquad\quad+\DhB{i}{h_i}\Big)\dy,}
\end{split}
\end{equation}
\begin{equation}\label{dF22}
\begin{split}
\partial_{r_i}\mathcal F_{i22}h_i=&\fit[.85\hsize]{{}-\kappa_i\ep b_i\fint\frac{\Big(h_i'(x)\big(r_i(x)-r_i(y)\big)+r_i'(x)\big(h_i(x)-h_i(y)\big)\Big)\cos(x-y)\big(1+\ep b_is_i(y)\big)}{\Dself{i}}\dy}\\
&\fit[.85\hsize]{{}+\kappa_i\ep b_i\fint\frac{r_i'(x)\big(r_i(x)-r_i(y)\big)\cos(x-y)\big(1+\ep b_is_i(y)\big)}{\big(\Dself{i}\big)^2}}\\
&\fit[.85\hsize]{\qquad\times\Big(\DhA{i}{h_i}}\\
&\fit[.85\hsize]{\qquad\quad+\DhB{i}{h_i}\Big)\dy,}
\end{split}
\end{equation}
\begin{equation}\label{dF23}
\begin{split}
\partial_{r_i}\mathcal F_{i23}h_i=&\fit[.85\hsize]{\ \kappa_i\fint\frac{\Big(\ep b_ih_i(x)\big(r_i(x)-r_i(y)\big)+\big(1+\ep b_ir_i(x)\big)\big(h_i(x)-h_i(y)\big)\Big)\sin(x-y)\big(1+\ep b_is_i(y)\big)}{\Dself{i}}\dy}\\
&\fit[.85\hsize]{{}-\kappa_i\fint\frac{\big(1+\ep b_ir_i(x)\big)\big(r_i(x)-r_i(y)\big)\sin(x-y)\big(1+\ep b_is_i(y)\big)}{\big(\Dself{i}\big)^2}}\\
&\fit[.85\hsize]{\qquad\times\Big(\DhA{i}{h_i}}\\
&\fit[.85\hsize]{\qquad\quad+\DhB{i}{h_i}\Big)\dy,}
\end{split}
\end{equation}
and
\begin{equation}\label{dF24}
\begin{split}
\partial_{r_i}\mathcal F_{i24}h_i=&\fit[.85\hsize]{{}-2\kappa_i\fint\frac{\big(1+\ep b_ir_i(x)\big)h_i(x)\sin(x-y)\big(1+\ep b_is_i(y)\big)}{\Dself{i}}\dy}\\
&\fit[.85\hsize]{{}+\kappa_i\fint\frac{\big(1+\ep b_ir_i(x)\big)^2\sin(x-y)\big(1+\ep b_is_i(y)\big)}{\big(\Dself{i}\big)^2}}\\
&\qquad\times\Big(2\ep b_i\big(r_i(x)-r_i(y)\big)\big(h_i(x)-h_i(y)\big)\\
&\qquad\quad+4\big(h_i(x)(1+\ep b_ir_i(y))+(1+\ep b_ir_i(x))h_i(y)\big)\sin^2\Big(\frac{x-y}{2}\Big)\Big)\dy .
\end{split}
\end{equation}
 Similarly, the derivative of the term $\mathcal F_{i2}$ with respect to $s_i$ is
\begin{equation}\label{dgF2}
\begin{split}
\partial_{s_i}\mathcal F_{i2}\eta_i=&\fit[.85\hsize]{{}-\kappa_i\ep b_i\fint\frac{\big(1+\ep b_ir_i(x)\big)r_i'(x)\big(1-\cos(x-y)\big)\eta_i(y)}{\Dself{i}}\dy}\\
&\fit[.85\hsize]{{}-\kappa_i\ep^2b_i^2\fint\frac{r_i'(x)\big(r_i(x)-r_i(y)\big)\cos(x-y)\eta_i(y)}{\Dself{i}}\dy}\\
&\fit[.85\hsize]{{}+\kappa_i\ep b_i\fint\frac{\big(1+\ep b_ir_i(x)\big)\big(r_i(x)-r_i(y)\big)\sin(x-y)\eta_i(y)}{\Dself{i}}\dy}\\
&\fit[.85\hsize]{{}-\kappa_i\fint\frac{\big(1+\ep b_ir_i(x)\big)^2\sin(x-y)\eta_i(y)}{\Dself{i}}\dy .}
\end{split}
\end{equation}
Now, we compute the Gateaux derivative of the functional \eqref{H2} with respect to $r_i$, which is given by
\begin{equation}\label{dH21}
\begin{split}
\partial_{r_i}\tilde{\mathcal H}_{i21}h_i=&\fit[.85\hsize]{\ 2\kappa_i\,\frac{\big(1+\ep b_is_i(x)\big)\Big(\big(1+\ep b_ir_i(x)\big)h_i(x)+\ep b_ir_i'(x)h_i'(x)\Big)}{\Big(\Ufull{i}\Big)^2}}\\
&\fit[.85\hsize]{\qquad\times\fint\frac{\big(1+\ep b_ir_i(x)\big)\big(1+\ep b_ir_i(y)\big)\big(\cos(x-y)-1\big)\big(1+\ep b_is_i(y)\big)}{\Dself{i}}\dy}\\
&\fit[.85\hsize]{{}-\kappa_i\,\frac{1+\ep b_is_i(x)}{\Ufull{i}}}\\
&\fit[.85\hsize]{\qquad\times\fint\frac{\Big(h_i(x)\big(1+\ep b_ir_i(y)\big)+\big(1+\ep b_ir_i(x)\big)h_i(y)\Big)\big(\cos(x-y)-1\big)\big(1+\ep b_is_i(y)\big)}{\Dself{i}}\dy}\\
&\fit[.85\hsize]{{}+\kappa_i\,\frac{1+\ep b_is_i(x)}{\Ufull{i}}}\\
&\fit[.85\hsize]{\qquad\times\fint\frac{\big(1+\ep b_ir_i(x)\big)\big(1+\ep b_ir_i(y)\big)\big(\cos(x-y)-1\big)\big(1+\ep b_is_i(y)\big)}{\big(\Dself{i}\big)^2}}\\
&\qquad\times\Big(2\ep b_i\big(r_i(x)-r_i(y)\big)\big(h_i(x)-h_i(y)\big)\\
&\qquad\quad+4\Big(h_i(x)\big(1+\ep b_ir_i(y)\big)+\big(1+\ep b_ir_i(x)\big)h_i(y)\Big)\sin^2\Big(\frac{x-y}{2}\Big)\Big)\dy .
\end{split}
\end{equation}
\begin{equation}\label{dH22}
\begin{split}
\partial_{r_i}\tilde{\mathcal H}_{i22}h_i=&\fit[.85\hsize]{\ \kappa_i\ep b_i\,\partial_{r_i}\Big(\frac{1+\ep b_is_i(x)}{\Ufull{i}}\Big)h_i}\\
&\fit[.85\hsize]{\qquad\times\fint\frac{r_i'(x)\big(r_i(x)-r_i(y)\big)\sin(x-y)\big(1+\ep b_is_i(y)\big)}{\Dself{i}}\dy}\\
&\fit[.85\hsize]{{}+\kappa_i\ep b_i\,\frac{1+\ep b_is_i(x)}{\Ufull{i}}}\\
&\fit[.85\hsize]{\qquad\times\fint\frac{\Big(h_i'(x)\big(r_i(x)-r_i(y)\big)+r_i'(x)\big(h_i(x)-h_i(y)\big)\Big)\sin(x-y)\big(1+\ep b_is_i(y)\big)}{\Dself{i}}\dy}\\
&\fit[.85\hsize]{{}-\kappa_i\ep b_i\,\frac{1+\ep b_is_i(x)}{\Ufull{i}}}\\
&\fit[.85\hsize]{\qquad\times\fint\frac{r_i'(x)\big(r_i(x)-r_i(y)\big)\sin(x-y)\big(1+\ep b_is_i(y)\big)}{\big(\Dself{i}\big)^2}}\\
&\fit[.85\hsize]{\qquad\times\Big(\DhA{i}{h_i}}\\
&\fit[.85\hsize]{\qquad\quad+\DhB{i}{h_i}\Big)\dy,}
\end{split}
\end{equation}
\begin{equation}\label{dH23}
\begin{split}
\partial_{r_i}\tilde{\mathcal H}_{i23}h_i=&\fit[.85\hsize]{\ \kappa_i\,\partial_{r_i}\Big(\frac{1+\ep b_is_i(x)}{\Ufull{i}}\Big)h_i}\\
&\fit[.85\hsize]{\qquad\times\fint\frac{\big(1+\ep b_ir_i(x)\big)\big(r_i(x)-r_i(y)\big)\big(1+\ep b_is_i(y)\big)}{\Dself{i}}\dy}\\
&\fit[.85\hsize]{{}+\kappa_i\,\frac{1+\ep b_is_i(x)}{\Ufull{i}}}\\
&\fit[.85\hsize]{\qquad\times\fint\frac{\Big(\ep b_ih_i(x)\big(r_i(x)-r_i(y)\big)+\big(1+\ep b_ir_i(x)\big)\big(h_i(x)-h_i(y)\big)\Big)\big(1+\ep b_is_i(y)\big)}{\Dself{i}}\dy}\\
&\fit[.85\hsize]{{}-\kappa_i\,\frac{1+\ep b_is_i(x)}{\Ufull{i}}}\\
&\fit[.85\hsize]{\qquad\times\fint\frac{\big(1+\ep b_ir_i(x)\big)\big(r_i(x)-r_i(y)\big)\big(1+\ep b_is_i(y)\big)}{\big(\Dself{i}\big)^2}}\\
&\fit[.85\hsize]{\qquad\times\Big(\DhA{i}{h_i}}\\
&\fit[.85\hsize]{\qquad\quad+\DhB{i}{h_i}\Big)\dy,}
\end{split}
\end{equation}
\begin{equation}\label{dH24}
\begin{split}
\partial_{r_i}\tilde{\mathcal H}_{i24}h_i=&\fit[.85\hsize]{{}-\kappa_i\bigg(\partial_{r_i}\Big(\frac{1+\ep b_is_i(x)}{\Ufull{i}}\Big)h_i\,r_i'(x)\big(1+\ep b_ir_i(x)\big)}\\
&\fit[.85\hsize]{\qquad+\frac{1+\ep b_is_i(x)}{\Ufull{i}}\Big(\ep b_ih_i(x)r_i'(x)+\big(1+\ep b_ir_i(x)\big)h_i'(x)\Big)\bigg)}\\
&\fit[.85\hsize]{\qquad\times\fint\frac{\sin(x-y)\big(1+\ep b_is_i(y)\big)}{\Dself{i}}\dy}\\
&\fit[.85\hsize]{{}+\kappa_i\,\frac{1+\ep b_is_i(x)}{\Ufull{i}}\,r_i'(x)\big(1+\ep b_ir_i(x)\big)}\\
&\fit[.85\hsize]{\qquad\times\fint\frac{\sin(x-y)\big(1+\ep b_is_i(y)\big)}{\big(\Dself{i}\big)^2}}\\
&\fit[.85\hsize]{\qquad\times\Big(\DhA{i}{h_i}}\\
&\fit[.85\hsize]{\qquad\quad+\DhB{i}{h_i}\Big)\dy .}
\end{split}
\end{equation}
Similarly, differentiating  the first term of \eqref{H2} with respect to $s_i$, we have
\begin{equation}\label{dgH21}
\begin{split}
\partial_{s_i}\tilde{\mathcal H}_{i21}\eta_i=&\fit[.85\hsize]{{}-\kappa_i\,\frac{\eta_i(x)}{\Ufull{i}}\fint\frac{\big(1+\ep b_ir_i(x)\big)\big(1+\ep b_ir_i(y)\big)\big(\cos(x-y)-1\big)\big(1+\ep b_is_i(y)\big)}{\Dself{i}}\dy}\\
&\fit[.85\hsize]{{}-\kappa_i\,\frac{1+\ep b_is_i(x)}{\Ufull{i}}\fint\frac{\big(1+\ep b_ir_i(x)\big)\big(1+\ep b_ir_i(y)\big)\big(\cos(x-y)-1\big)\eta_i(y)}{\Dself{i}}\dy. }
\end{split}
\end{equation}

 By similar computations, the derivatives of the interaction terms are obtained by differentiating the integrands in \eqref{F3} and \eqref{H3}. Hence
\begin{equation*}
\partial_{r_{3-i}}\mathcal F_{i3}h_{3-i}=\ep\kappa_{3-i}\mathcal \partial_{r_{3-i}}\mathcal Q_{i3} .
\end{equation*}
and 
\begin{equation*}
\partial_{r_{3-i}}\mathcal F_{i4}h_{3-i}=\ep\kappa_{3-i}\mathcal \partial_{r_{3-i}}\mathcal Q_{i4} .
\end{equation*}

Now, we prove that the Gateaux derivatives are well defined. To this aim, we show
\begin{equation}\label{gateaux}
\fit[\hsize]{\lim_{t\to0}\norm{\frac{\mathcal F_{i2j}(\ep,\Omega,x_0,\ldots,r_i+th_i,\ldots)-\mathcal F_{i2j}(\ep,\Omega,x_0,r_1,s_1,r_2,s_2)}{t}-\partial_{r_i}\mathcal F_{i2j}h_i}_{Y^k}=0,\qquad j=1,2,3,4,}
\end{equation}
and similarly for the other terms. Notice  the most singular term is $\mathcal F_{i22}$. We use the notation introduced in the proof of the continuity in Lemma \ref{lem:cont}, and we denote by
\begin{equation*}
\fit[\hsize]{D_\ep'(r_i)[h_i]:=2\ep^2b_i^2\,\Delta r_i\,\Delta h_i+4\ep b_i\big(h_i(1+\ep b_i\tilde r_i)+(1+\ep b_ir_i)\tilde h_i\big)\sin^2\Big(\frac{x-y}{2}\Big)}
\end{equation*}
the derivative of $D_\ep(r_i)$ in the direction $h_i$.  Then, we have
\begin{equation*}
\begin{split}
&\frac{\mathcal F_{i22}(\ep,r_i+th_i,s_i)-\mathcal F_{i22}(\ep,r_i,s_i)}{t}-\partial_{r_i}\mathcal F_{i22}(\ep,r_i,s_i)h_i\\
&\fit[.85\hsize]{{}=-\kappa_i\ep b_i\fint r_i'\,\Delta r_i\cos(x-y)\big(1+\ep b_i\tilde s_i\big)\Big(\frac{1}{t}\Big(\frac{1}{D_\ep(r_i+th_i)}-\frac{1}{D_\ep(r_i)}\Big)+\frac{D_\ep'(r_i)[h_i]}{D_\ep(r_i)^2}\Big)\dy}\\
&\fit[.85\hsize]{\quad-\kappa_i\ep b_i\fint\big(h_i'\,\Delta r_i+r_i'\,\Delta h_i\big)\cos(x-y)\big(1+\ep b_i\tilde s_i\big)\Big(\frac{1}{D_\ep(r_i+th_i)}-\frac{1}{D_\ep(r_i)}\Big)\dy}\\
&\fit[.85\hsize]{\quad-\kappa_i\ep b_it\fint\frac{h_i'\,\Delta h_i\cos(x-y)\big(1+\ep b_i\tilde s_i\big)}{D_\ep(r_i+th_i)}\dy=:G_{i21}+G_{i22}+G_{i23} .}
\end{split}
\end{equation*}
Based on \eqref{diffD}, the terms $G_{i21}$ and $G_{i23}$ satisfy the bounds
    \begin{equation*}
    	\|\partial^kG_{i21}(x\pm a\mathbf{i})\|_{L^2}\le tC\|h_i\|_{X^{k+1}}\|r_i\|_{X^{k+1}}(1+\|s_i\|_{X^k})(1+\|r_i\|_{X^{k+1}})\|h_i\|_{X^k},
    \end{equation*}
and
\begin{equation*}
    	\|\partial^kG_{i23}(x\pm a\mathbf{i})\|_{L^2}\le tC\|r_i\|_{X^{k+1}}\|r_i\|_{X^{k+1}}(1+\|s_i\|_{X^k})(1+\|r_i\|_{X^{k+1}})\|h_i\|_{X^k}.
    \end{equation*}
For $G_{i22}$, an application of the mean value theorem gives
 \begin{equation*}
\frac{1}{D_\ep(r_i+th_i)}-\frac{1}{D_\ep(r_i)}+t\frac{D_\ep'(r_i)[h_i]}{D_\ep(r_i)^2}\le \frac{t^2}{4\sin^2\left(\frac{x-y}{2}\right)} \mathsf b(\varepsilon,r_i,h_i),
    \end{equation*}
with $\Vert{}\mathsf b(\varepsilon,r_i,h_i)(x\pm a\mathbf{i})\Vert{}_{L^\infty} \le C(\varepsilon,\Vert{}r_i\Vert{}_{X^{k+1}},\Vert{}h_i\Vert{}_{x^{k+1}}) < \infty$. This implies that
    \begin{equation*}
    	\|\partial^kG_{i22}(x\pm a\mathbf{i})\|_{L^2}\le tC\|r_i\|_{X^{k+1}}\|r_i\|_{X^{k+1}}(1+\|s_i\|_{X^k})\|\xi(\varepsilon,r_i,h_i)(x\pm a\mathbf{i})\|_{L^\infty}.
    \end{equation*}
Regarding $G_{i22}$, applying the mean value theorem yields
    \begin{equation*}
    	\|\partial^kG_{i24}(x\pm a\mathbf{i})\|_{L^2}\le tC\|r_i\|_{X^{k+1}}\|h_i\|_{X^{k+1}}(1+\|s_i\|_{X^k}).
    \end{equation*}
%%%%5
 Finally, for the term $G_{i24}$, it readily follows that
\begin{equation*}
\frac{1}{D_\ep(r_i+th_i)}-\frac{1}{D_\ep(r_i)}+t\frac{D_\ep'(r_i)[h_i]}{D_\ep(r_i)^2}\le \frac{t^2}{4\sin^2\left(\frac{x-y}{2}\right)} \mathsf b(\varepsilon,r_i,h_i),
    \end{equation*}
    and
    \begin{equation*}
    	\|\partial^kG_{i23}(x\pm a\mathbf{i})\|_{L^2}\le tC\|r_i\|_{X^{k+1}}\|r_i\|_{X^{k+1}}(1+\|s_i\|_{X^k})(1+\|r_i\|_{X^{k+1}})\|h_i\|_{X^k}.
    \end{equation*}
    For $G_{i22}$, we use mean value theorem, to get
    \begin{equation*}
\frac{1}{D_\ep(r_i+th_i)}-\frac{1}{D_\ep(r_i)}+t\frac{D_\ep'(r_i)[h_i]}{D_\ep(r_i)^2}\le \frac{t^2}{4\sin^2\left(\frac{x-y}{2}\right)} \mathsf b(\varepsilon,r_i,h_i),
    \end{equation*}
    where $\|\mathsf b(\varepsilon,r_i,h_i)(x\pm a\mathbf{i})\|_{L^\infty}\le C(\varepsilon,\|r_i\|_{X^{k+1}},\|h_i\|_{X^{k+1}})<\infty$. As a result, it holds
    \begin{equation*}
    	\|\partial^kG_{i22}(x\pm a\mathbf{i})\|_{L^2}\le tC\|r_i\|_{X^{k+1}}\|r_i\|_{X^{k+1}}(1+\|s_i\|_{X^k})\|\xi(\varepsilon,r_i,h_i)(x\pm a\mathbf{i})\|_{L^\infty}.
    \end{equation*}
    For the  term $G_{i24}$, it is obvious that
    \begin{equation*}
    	\|\partial^kG_{i24}(x\pm a\mathbf{i})\|_{L^2}\le tC\|r_i\|_{X^{k+1}}\|h_i\|_{X^{k+1}}(1+\|s_i\|_{X^k}).
    \end{equation*}
   Then, by sending $t\to0$, we arrive at \eqref{gateaux} for $j=2$.

Notice that $\mathcal F_{i2}$ is affine in $s_i$, and $\tilde{\mathcal H}_{i2}$ is a polynomial of degree two in $s_i$. The derivatives with respect to $s_i$. Now we prove the continuity of the derivatives. Let $h=(\beta_1,\beta_2,h_1,\eta_1,h_2,\eta_2)\in\mathcal X^k$ be a fixed direction. The expression $D_g\mathcal \mathcal F(\ep,g)h$ is a sum of the integrals \eqref{dOmega}--\eqref{dgH21}, which are linear in $h$, and whose coefficients depend on $(\ep,g)$ through the same quantities as in Lemma \ref{lem:cont}. 
The continuity of $\partial_{r_i} \mathcal F_{i2}h_i$, as well as that of $\partial_{r_i}\mathcal F_{i3}$, $\partial_{r_i}\mathcal F_{i4}$, follows by the same arguments established in Lemma \ref{lem:cont}. Therefore, we omit the details. Repeating the arguments of the proof of Lemma \ref{lem:cont}, we can conclude that
$ D_g\mathcal \mathcal F$ is continuous. 

Concerning the derivative with respect to $\ep$, as shown in the proof of the expansions in Lemma \ref{lem:cont}, the partial derivative $\partial_\ep \mathcal F$ exists and is given by integrals of the same type as those studied in Lemma \ref{lem:cont}. Arguing as in the proof of the continuity in that lemma, it is continuous in $(\ep,g)$. Together with the continuity of $D_g\mathcal F$, this shows that $\mathcal F$ is of class $C^1$.

%%%
%%%
%%
%%%
%%%
%%
It remains to compute the derivatives at $\ep=0$.  From \eqref{dF21}--\eqref{dgF2}, we obtain
\begin{equation*}
\begin{split}
\partial_{r_i}\mathcal F_{i2}(0,g)h_i&\fit[.85\hsize]{{}=-\kappa_i\fint\frac{h_i'(x)\,2\sin^2(\frac{x-y}{2})}{4\sin^2(\frac{x-y}{2})}\dy+\kappa_i\fint\frac{(h_i(x)-h_i(y))\sin(x-y)}{4\sin^2(\frac{x-y}{2})}\dy}\\
&\fit[.85\hsize]{\quad-2\kappa_i\fint\frac{h_i(x)\sin(x-y)}{4\sin^2(\frac{x-y}{2})}\dy+\kappa_i\fint\frac{\sin(x-y)\,4\big(h_i(x)+h_i(y)\big)\sin^2(\frac{x-y}{2})}{16\sin^4(\frac{x-y}{2})}\dy}\\
&\fit[.85\hsize]{{}=-\frac{\kappa_i}{2}h_i'(x)-\kappa_i\fint\frac{h_i(y)\sin(x-y)}{4\sin^2(\frac{x-y}{2})}\dy+0+\kappa_i\fint\frac{h_i(y)\sin(x-y)}{4\sin^2(\frac{x-y}{2})}\dy}\\
&=-\frac{\kappa_i}{2}h_i'(x).
\end{split}
\end{equation*}
Then, we have
\begin{equation}
    \partial_{r_i}\mathcal F_i(0,g)h_i=-\frac{\kappa_i}{2}h_i'(x)
\end{equation}
Similary, we have 
\begin{equation}
\partial_{s_i}\mathcal F_{i}(0,g)\eta_i =\kappa_i\fint\frac{\eta_i(y)\sin(x-y)}{4\sin^2(\frac{x-y}{2})}\dy . 
\end{equation}
Now, applying \eqref{dH21} yields
\begin{equation}
\partial_{r_i}\mathcal H_{i21}(0,g)h_i=\frac{\kappa_i}{2}(I-\Pm)\big[-h_i-\fint(h_i(x)+h_i(y))\dy\big]=-\kappa_ih_i .
\end{equation}
In view of \eqref{dH22}--\eqref{dH24}, the only other non-vanishing term at $\ep=0$ is given by
\begin{equation*}
\fit[\hsize]{\partial_{r_i}\tilde{\mathcal H}_{i23}(0,g)h_i=\kappa_i\fint\frac{h_i(x)-h_i(y)}{4\sin^2(\frac{x-y}{2})}\dy,}
\end{equation*}
since $\fint\frac{\sin(x-y)}{4\sin^2((x-y)/2)}\dy=0$. Consequently, we have
\begin{equation*}
\fit[\hsize]{\partial_{r_i}\mathcal H_{i2}(0,g)h_i=-\kappa_ih_i(x)+\kappa_i\fint\frac{h_i(x)-h_i(y)}{4\sin^2(\frac{x-y}{2})}\dy,\qquad \partial_{s_i}\mathcal H_{i2}(0,g)\eta_i=\frac{\kappa_i}{2}\eta_i(x),}
\end{equation*}
where the second identity derives from \eqref{dgH21} evaluated at $\ep=0$. Then, one gets
\begin{equation}
\partial_{r_i}\mathcal H_i(0,g)h_i=-\kappa_ih_i(x)+\kappa_i\fint\frac{h_i(x)-h_i(y)}{4\sin^2(\frac{x-y}{2})}\dy .
\end{equation}
Analogous computations lead to
\begin{equation}
\partial_{s_i}\mathcal H_i(0,g)\eta_i=\frac{\kappa_i}{2}\eta_i(x).
\end{equation}
In the same way we have
\begin{equation}
\partial_{r_{3-i}}\mathcal F_i(0,g)=\partial_{s_{3-i}}\mathcal F_i(0,g)=\partial_{r_{3-i}}\mathcal H_i(0,g)=\partial_{s_{3-i}}\mathcal H_i(0,g)=0.
\end{equation}
Finally, evaluating \eqref{dOmega} at $\ep=0$ provides the derivatives with respect to $\Omega$ and $x_0$
\begin{equation}
    \begin{split}
\partial_\Omega \mathcal F_i(0,g)=x_{0,i}\sin x,\qquad \partial_{x_0}\mathcal F_i(0,g)=(-1)^{i+1}\Omega\sin x,\\
\partial_\Omega \mathcal H_i(0,g)=x_{0,i}\cos x,\qquad \partial_{x_0}\mathcal H_i(0,g)=(-1)^{i+1}\Omega\cos x .
    \end{split}
\end{equation}
This proves \eqref{lin0} and completes the proof of Lemma \ref{lem:C1}.
\end{proof}

\subsection{Point vortex pairs and the linearized operator}

\begin{lemma}\label{critical}
For every $(\Omega,x_0)\in\R^2$ and for $i=1,2$ we have, 
\begin{equation}\label{Fat0}
\fit[\hsize]{\mathcal F_i(0,\Omega,x_0,0,0,0,0)=\Big(\Omega x_{0,i}-\frac{\kappa_{3-i}}{d}\Big)\sin x,\qquad \mathcal H_i(0,\Omega,x_0,0,0,0,0)=\Big(\Omega x_{0,i}-\frac{\kappa_{3-i}}{d}\Big)\cos x .}
\end{equation}
Consequently, if $\kappa_1+\kappa_2\neq0$, then $\bfF(0,\Omega,x_0,0,0,0,0)=0$ if and only if $(\Omega,x_0)=(\Omega^*,x_0^*)$, where
\begin{equation}\label{vel}
\Omega^*=\frac{\kappa_1+\kappa_2}{d^2},\qquad x_0^*=\frac{d\kappa_2}{\kappa_1+\kappa_2}.
\end{equation}
In other words, the point vortices $2\pi\kappa_1\boldsymbol\delta_0$ and $2\pi\kappa_2\boldsymbol\delta_{d\bfe_1}$ rotate uniformly about the centroid $(x_0^*,0)$ with angular velocity $\Omega^*$.
\end{lemma}

\begin{proof}
Identity \eqref{Fat0} follows directly from \eqref{Fexp}--\eqref{Hexp} by setting $\ep=0$ and $r_i=s_i=0$, with $x_{0,i}$ given by \eqref{xbi}. Therefore,
$$\bfF(0,\Omega,x_0,0,0,0,0)=\big(\mathcal F_1,\mathcal H_1,\mathcal F_2,\mathcal H_2\big)(0,\Omega,x_0,0,0,0,0)=0. $$
if and only if
\begin{equation}\label{f1f2}
\bfF_1(0,\Omega,x_0,0,0,0,0)=\left(\Omega x_0-\frac{\kappa_2}{d}\right)\sin (x)=0,
\end{equation}
\begin{equation}\label{f1f2b}
\bfF_2(0,\Omega,x_0,0,0,0,0)=\left(\Omega(d-x_0)-\frac{\kappa_1}{d}\right)\sin (x)=0, 
\end{equation}
\begin{equation}\label{f1f2c}
\mathcal H_1(0,\Omega,x_0,0,0,0,0)=\left(\Omega x_0-\frac{\kappa_2}{d}\right)\cos (x)=0,
\end{equation}
\begin{equation}\label{f1f2d}
\mathcal  H_2(0,\Omega,x_0,0,0,0,0)=\left(\Omega(d-x_0)-\frac{\kappa_1}{d}\right)\cos (x)=0, 
\end{equation}
Since $\sin(x)$ and $\cos(x)$ do not vanish simultaneously, equations \eqref{f1f2}--\eqref{f1f2d} imply that their respective coefficients must be identically zero. This yields the algebraic system
\begin{equation}\label{eq:system_Omega_x0}
\Omega x_0 = \frac{\kappa_2}{d} \quad \text{and} \quad \Omega(d-x_0) = \frac{\kappa_1}{d}.
\end{equation}
Adding both equations in \eqref{eq:system_Omega_x0} allows us to eliminate $x_0$, leading to
\begin{equation*}
\Omega = \frac{\kappa_1+\kappa_2}{d^2} := \Omega^*.
\end{equation*}
Assuming $\kappa_1+\kappa_2 \neq 0$ so that $\Omega^* \neq 0$, we can substitute this value  into the first equation of \eqref{eq:system_Omega_x0} to obtain
\begin{equation}
x_0 = \frac{\kappa_2 d}{\kappa_1+\kappa_2} := x_0^*.
\end{equation}
\end{proof}

Take $\hfr=(\beta_1,\beta_2,h_1,\eta_1,h_2,\eta_2)\in\mathcal X^k$, where
\begin{equation}\label{hexp1}
h_i(x)=\sum_{j=1}^\infty a_j^i\cos(jx),\qquad \eta_i(x)=\sum_{j=1}^\infty b_j^i\cos(jx),\qquad b_1^i=-a_1^i,\quad i=1,2 .
\end{equation}
Then the linearization of $\bfF$ at $(0,\gfr_0)$ has the following Fourier series form
\begin{equation}\label{DF1}
\begin{split}
D_\gfr\bfF(0,\gfr_0)\hfr
&=\begin{pmatrix}
\partial_\Omega\mathcal F_1\,\beta_1+\partial_{x_0}\mathcal F_1\,\beta_2+\partial_{r_1}\mathcal F_1\,h_1+\partial_{s_1}\mathcal F_1\,\eta_1\\[2pt]
\partial_\Omega\mathcal H_1\,\beta_1+\partial_{x_0}\mathcal H_1\,\beta_2+\partial_{r_1}\mathcal H_1\,h_1+\partial_{s_1}\mathcal H_1\,\eta_1\\[2pt]
\partial_\Omega\mathcal F_2\,\beta_1+\partial_{x_0}\mathcal F_2\,\beta_2+\partial_{r_2}\mathcal F_2\,h_2+\partial_{s_2}\mathcal F_2\,\eta_2\\[2pt]
\partial_\Omega\mathcal H_2\,\beta_1+\partial_{x_0}\mathcal H_2\,\beta_2+\partial_{r_2}\mathcal H_2\,h_2+\partial_{s_2}\mathcal H_2\,\eta_2
\end{pmatrix}\\
&=\begin{pmatrix}\lambda_1\sin x\\ \lambda_1\cos x\\ \lambda_2\sin x\\ \lambda_2\cos x\end{pmatrix}
+\sum_{j=1}^\infty\begin{pmatrix}\hat A_j^1\sin(jx)\\ \hat B_j^1\cos(jx)\\ \hat A_j^2\sin(jx)\\ \hat B_j^2\cos(jx)\end{pmatrix} .
\end{split}
\end{equation}
Recall that all the partial derivatives are evaluated at $(0,\gfr_0)$ with $\lambda_i:=\beta_1x_{0,i}^*+(-1)^{i+1}\beta_2\Omega^*$  and
\begin{equation*}
\begin{pmatrix}\hat A_j^i\\ \hat B_j^i\end{pmatrix}=M_j^{(i)}\begin{pmatrix}a_j^i\\ b_j^i\end{pmatrix},\qquad i=1,2,
\end{equation*}
where $M_j^{(i)}$ a $2\times2$ matrix and $a_j^i,b_j^i$ given in \eqref{hexp}. In the following lemma we compute the matrices $M_j^{(i)}$ and prove that $D_\gfr\bfF(0,\gfr_0)$ is an isomorphism from $\mathcal X^k$ onto $\mathcal Y^k$. 

%%%%
%%%%
%%%%

%%%%
%%%%
%%%%
Take $\hfr=(\beta_1,\beta_2,h_1,\eta_1,h_2,\eta_2)\in\mathcal X^k$, where
\begin{equation}\label{hexp}
h_i(x)=\sum_{j=1}^\infty a_j^i\cos(jx),\qquad \eta_i(x)=\sum_{j=1}^\infty b_j^i\cos(jx),\qquad b_1^i=-a_1^i,\quad i=1,2 .
\end{equation}
By \eqref{lin0}, the derivatives of $\mathcal F_i$ and $\mathcal H_i$ with respect to $(r_{3-i},s_{3-i})$ vanish at $\ep=0$. Then the linearization of $\bfF$ at $(0,\gfr_0)$ has the following matrix form and Fourier series form
\begin{equation}\label{DF}
\begin{split}
D_\gfr\bfF(0,\gfr_0)\hfr
&=\begin{pmatrix}
\partial_\Omega\mathcal F_1 & \partial_{x_0}\mathcal F_1 & \partial_{r_1}\mathcal F_1 & \partial_{s_1}\mathcal F_1 & 0 & 0\\[2pt]
\partial_\Omega\mathcal H_1 & \partial_{x_0}\mathcal H_1 & \partial_{r_1}\mathcal H_1 & \partial_{s_1}\mathcal H_1 & 0 & 0\\[2pt]
\partial_\Omega\mathcal F_2 & \partial_{x_0}\mathcal F_2 & 0 & 0 & \partial_{r_2}\mathcal F_2 & \partial_{s_2}\mathcal F_2\\[2pt]
\partial_\Omega\mathcal H_2 & \partial_{x_0}\mathcal H_2 & 0 & 0 & \partial_{r_2}\mathcal H_2 & \partial_{s_2}\mathcal H_2
\end{pmatrix}
\begin{pmatrix}\beta_1\\ \beta_2\\ h_1\\ \eta_1\\ h_2\\ \eta_2\end{pmatrix}\\
&=\begin{pmatrix}\lambda_1\sin x\\ \lambda_1\cos x\\ \lambda_2\sin x\\ \lambda_2\cos x\end{pmatrix}
+\sum_{j=1}^\infty\begin{pmatrix}\hat A_j^1\sin(jx)\\ \hat B_j^1\cos(jx)\\ \hat A_j^2\sin(jx)\\ \hat B_j^2\cos(jx)\end{pmatrix},
\end{split}
\end{equation}
where all the partial derivatives are evaluated at $(0,\gfr_0)$,
\begin{equation}\label{lambdai}
\lambda_i:=\beta_1x_{0,i}^*+(-1)^{i+1}\beta_2\Omega^*,\qquad x_{0,1}^*:=x_0^*,\quad x_{0,2}^*:=d-x_0^*,
\end{equation}
and
\begin{equation*}
\begin{pmatrix}\hat A_j^i\\ \hat B_j^i\end{pmatrix}=M_j^{(i)}\begin{pmatrix}a_j^i\\ b_j^i\end{pmatrix},\qquad i=1,2,
\end{equation*}
with $M_j^{(i)}$ a $2\times2$ matrix and $a_j^i,b_j^i$ given in \eqref{hexp}. In the following lemma we compute the matrices $M_j^{(i)}$, we show that, in terms of the Fourier coefficients, $D_\gfr\bfF(0,\gfr_0)$ is an infinite block diagonal matrix with $4\times4$ blocks, and we prove that $D_\gfr\bfF(0,\gfr_0)$ is an isomorphism from $\mathcal X^k$ onto $\mathcal Y^k$.
\begin{lemma}\label{lem:iso}
Assume $\kappa_1+\kappa_2\ne0$, let $\hfr\in\mathcal X^k$ be as in \eqref{hexp}, and write
\begin{equation}\label{uvexp}
D_\gfr\bfF(0,\gfr_0)\hfr=(u_1,v_1,u_2,v_2),\qquad u_i(x)=\sum_{j=1}^\infty A_j^i\sin(jx),\qquad v_i(x)=\sum_{j=1}^\infty B_j^i\cos(jx).
\end{equation}
Then the following assertions hold true.
\begin{enumerate}
\item[(a)] For $i=1,2$ and $j\ge1$, the matrices in \eqref{DF} are
\begin{equation}\label{Mj}
M_j^{(i)}=\frac{\kappa_i}{2}\begin{pmatrix} j & -1\\ j-2 & 1\end{pmatrix}.
\end{equation}
\item[(b)] In terms of the Fourier coefficients, $D_\gfr\bfF(0,\gfr_0)$ is block diagonal. For every $j\ge1$, the Fourier coefficients of order $j$ of $D_\gfr\bfF(0,\gfr_0)\hfr$ depend only on the Fourier coefficients of order $j$ of $\hfr$, namely
\begin{equation}\label{DFblocks}
\begin{pmatrix}A_1^1\\ B_1^1\\ A_1^2\\ B_1^2\end{pmatrix}=\mathbb M_1\begin{pmatrix}\beta_1\\ \beta_2\\ a_1^1\\ a_1^2\end{pmatrix},\qquad
\begin{pmatrix}A_j^1\\ B_j^1\\ A_j^2\\ B_j^2\end{pmatrix}=\mathbb M_j\begin{pmatrix}a_j^1\\ b_j^1\\ a_j^2\\ b_j^2\end{pmatrix},\quad j\ge2,
\end{equation}
where
\begin{equation}\label{MM1}
\mathbb M_1:=\begin{pmatrix}
x_0^* & \Omega^* & \kappa_1 & 0\\
x_0^* & \Omega^* & -\kappa_1 & 0\\
d-x_0^* & -\Omega^* & 0 & \kappa_2\\
d-x_0^* & -\Omega^* & 0 & -\kappa_2
\end{pmatrix},\qquad
\mathbb M_j:=\begin{pmatrix}M_j^{(1)} & 0\\ 0 & M_j^{(2)}\end{pmatrix},\quad j\ge2 .
\end{equation}
Moreover,
\begin{equation}\label{detMM}
\det\mathbb M_1=\frac{4\kappa_1\kappa_2(\kappa_1+\kappa_2)}{d}\neq0,\qquad \det\mathbb M_j=\frac{\kappa_1^2\kappa_2^2(j-1)^2}{4}\neq0,\quad j\ge2 .
\end{equation}
\item[(c)] $D_\gfr\bfF(0,\gfr_0):\mathcal X^k\to\mathcal Y^k$ is an isomorphism. More precisely, given $(u_1,v_1,u_2,v_2)\in\mathcal Y^k$ as in \eqref{uvexp}, the unique $\hfr\in\mathcal X^k$ such that $D_\gfr\bfF(0,\gfr_0)\hfr=(u_1,v_1,u_2,v_2)$ is given by \eqref{hexp}, where $b_1^i=-a_1^i$,
\begin{equation}\label{DFinvblocks}
\begin{pmatrix}\beta_1\\ \beta_2\\ a_1^1\\ a_1^2\end{pmatrix}=\mathbb M_1^{-1}\begin{pmatrix}A_1^1\\ B_1^1\\ A_1^2\\ B_1^2\end{pmatrix},\qquad
\begin{pmatrix}a_j^1\\ b_j^1\\ a_j^2\\ b_j^2\end{pmatrix}=\mathbb M_j^{-1}\begin{pmatrix}A_j^1\\ B_j^1\\ A_j^2\\ B_j^2\end{pmatrix},\quad j\ge2,
\end{equation}
\begin{equation}\label{MM1inv}
\mathbb M_1^{-1}=\frac{1}{2d}\begin{pmatrix}
1 & 1 & 1 & 1\\[2pt]
\frac{d-x_0^*}{\Omega^*} & \frac{d-x_0^*}{\Omega^*} & -\frac{x_0^*}{\Omega^*} & -\frac{x_0^*}{\Omega^*}\\[4pt]
\frac{d}{\kappa_1} & -\frac{d}{\kappa_1} & 0 & 0\\[4pt]
0 & 0 & \frac{d}{\kappa_2} & -\frac{d}{\kappa_2}
\end{pmatrix},\qquad
\mathbb M_j^{-1}=\begin{pmatrix}\big(M_j^{(1)}\big)^{-1} & 0\\ 0 & \big(M_j^{(2)}\big)^{-1}\end{pmatrix},\quad j\ge2,
\end{equation}
and
\begin{equation}\label{Mjinv}
\big(M_j^{(i)}\big)^{-1}=\frac{1}{\kappa_i(j-1)}\begin{pmatrix}1 & 1\\ 2-j & j\end{pmatrix},\qquad j\ge2,\quad i=1,2 .
\end{equation}
\end{enumerate}
\end{lemma}

\begin{proof}
By Lemma \ref{lem:C1}, the entries of the matrix in \eqref{DF} are given by \eqref{lin0} with $(\Omega,x_0)=(\Omega^*,x_0^*)$, namely, for $i=1,2$,
\begin{equation*}
\begin{gathered}
\partial_\Omega\mathcal F_i(0,\gfr_0)=x_{0,i}^*\sin x,\qquad \partial_{x_0}\mathcal F_i(0,\gfr_0)=(-1)^{i+1}\Omega^*\sin x,\\
\partial_\Omega\mathcal H_i(0,\gfr_0)=x_{0,i}^*\cos x,\qquad \partial_{x_0}\mathcal H_i(0,\gfr_0)=(-1)^{i+1}\Omega^*\cos x,\\
\fit[\hsize]{\partial_{r_i}\mathcal F_i(0,\gfr_0)h_i=-\frac{\kappa_i}{2}h_i'(x),\qquad \partial_{s_i}\mathcal F_i(0,\gfr_0)\eta_i=-\kappa_i\fint\frac{\eta_i(y)\sin(x-y)}{4\sin^2(\frac{x-y}{2})}\dy,}\\
\fit[\hsize]{\partial_{r_i}\mathcal H_i(0,\gfr_0)h_i=-\kappa_ih_i(x)+\kappa_i\fint\frac{h_i(x)-h_i(y)}{4\sin^2(\frac{x-y}{2})}\dy,\qquad \partial_{s_i}\mathcal H_i(0,\gfr_0)\eta_i=\frac{\kappa_i}{2}\eta_i(x).}
\end{gathered}
\end{equation*}
Hence the contribution of the first two columns of the matrix in \eqref{DF} is $\lambda_i\sin x$ in the rows of $\mathcal F_i$ and $\lambda_i\cos x$ in the rows of $\mathcal H_i$, with $\lambda_i$ given in \eqref{lambdai}. Moreover, by \eqref{hexp} and \eqref{HLcos}, $h_i'(x)=-\sum_{j\ge1}ja_j^i\sin(jx)$,
\begin{equation*}
\fint\frac{\eta_i(y)\sin(x-y)}{4\sin^2(\frac{x-y}{2})}\dy=\frac12\sum_{j\ge1}b_j^i\sin(jx),\qquad \fint\frac{h_i(x)-h_i(y)}{4\sin^2(\frac{x-y}{2})}\dy=\frac12\sum_{j\ge1}ja_j^i\cos(jx),
\end{equation*}
and therefore
\begin{equation*}
\begin{split}
\partial_{r_i}\mathcal F_i(0,\gfr_0)h_i+\partial_{s_i}\mathcal F_i(0,\gfr_0)\eta_i&=\frac{\kappa_i}{2}\sum_{j=1}^\infty\big(ja_j^i-b_j^i\big)\sin(jx),\\
\partial_{r_i}\mathcal H_i(0,\gfr_0)h_i+\partial_{s_i}\mathcal H_i(0,\gfr_0)\eta_i&=\frac{\kappa_i}{2}\sum_{j=1}^\infty\big((j-2)a_j^i+b_j^i\big)\cos(jx).
\end{split}
\end{equation*}
That is, $\hat A_j^i=\frac{\kappa_i}{2}\big(ja_j^i-b_j^i\big)$ and $\hat B_j^i=\frac{\kappa_i}{2}\big((j-2)a_j^i+b_j^i\big)$, which proves \eqref{Mj}.

Comparing \eqref{DF} with \eqref{uvexp}, we get, for $i=1,2$,
\begin{equation*}
A_1^i=\lambda_i+\hat A_1^i,\qquad B_1^i=\lambda_i+\hat B_1^i,\qquad A_j^i=\hat A_j^i,\quad B_j^i=\hat B_j^i,\quad j\ge2 .
\end{equation*}
Since $b_1^i=-a_1^i$, \eqref{Mj} gives $\hat A_1^i=\frac{\kappa_i}{2}\big(a_1^i-b_1^i\big)=\kappa_ia_1^i$ and $\hat B_1^i=\frac{\kappa_i}{2}\big(-a_1^i+b_1^i\big)=-\kappa_ia_1^i$, so that
\begin{equation*}
A_1^i=\lambda_i+\kappa_ia_1^i,\qquad B_1^i=\lambda_i-\kappa_ia_1^i,\qquad i=1,2 .
\end{equation*}
By \eqref{lambdai}, this is the first identity in \eqref{DFblocks}, and the second one follows from $A_j^i=\hat A_j^i$, $B_j^i=\hat B_j^i$ and \eqref{Mj}.

 Notice that the matrix
\begin{equation*}
M_1^{(i)}=\frac{\kappa_i}{2}\begin{pmatrix}1&-1\\-1&1\end{pmatrix}
\end{equation*}
is not invertible. Its kernel is spanned by $(1,1)$, which corresponds to the infinitesimal translations of the $i$-th sheet in the direction $\bfe_1$, and its range is spanned by $(1,-1)$. This is the reason why we restrict the first Fourier coefficients by the condition $b_1^i=-a_1^i$ in \eqref{Ak} and add the unknowns $\Omega$ and $x_0$. Indeed, the columns of $\mathbb M_1$ corresponding to $a_1^1$ and $a_1^2$ are $\kappa_1(1,-1,0,0)^T$ and $\kappa_2(0,0,1,-1)^T$, that is, they span the ranges of $M_1^{(1)}$ and $M_1^{(2)}$. The columns corresponding to $\beta_1$ and $\beta_2$ are
\begin{equation*}
x_0^*(1,1,0,0)^T+(d-x_0^*)(0,0,1,1)^T\qquad\text{and}\qquad \Omega^*\big((1,1,0,0)^T-(0,0,1,1)^T\big),
\end{equation*}
and they span the two complementary directions $(1,1,0,0)^T$ and $(0,0,1,1)^T$ if and only if
\begin{equation*}
\det\begin{pmatrix}x_0^* & \Omega^*\\ d-x_0^* & -\Omega^*\end{pmatrix}=-x_0^*\Omega^*-(d-x_0^*)\Omega^*=-d\Omega^*=-\frac{\kappa_1+\kappa_2}{d}\neq0,
\end{equation*}
that is, if and only if $\kappa_1+\kappa_2\neq0$.

 Subtracting the second row of $\mathbb M_1$ from the first one and the fourth row from the third one, which does not change the determinant, we obtain
\begin{equation*}
\det\mathbb M_1=\det\begin{pmatrix}
0 & 0 & 2\kappa_1 & 0\\
x_0^* & \Omega^* & -\kappa_1 & 0\\
0 & 0 & 0 & 2\kappa_2\\
d-x_0^* & -\Omega^* & 0 & -\kappa_2
\end{pmatrix}.
\end{equation*}
Expanding along the first row and then along the row containing $2\kappa_2$, we get
\begin{equation*}
\det\mathbb M_1=2\kappa_1\cdot(-2\kappa_2)\det\begin{pmatrix}x_0^* & \Omega^*\\ d-x_0^* & -\Omega^*\end{pmatrix}=4\kappa_1\kappa_2\,d\,\Omega^*=\frac{4\kappa_1\kappa_2(\kappa_1+\kappa_2)}{d},
\end{equation*}
by \eqref{vel}. Moreover, $\det M_j^{(i)}=\frac{\kappa_i^2}{4}\big(j+(j-2)\big)=\frac{\kappa_i^2}{2}(j-1)$, so that $\det\mathbb M_j=\det M_j^{(1)}\det M_j^{(2)}=\frac{\kappa_1^2\kappa_2^2(j-1)^2}{4}$. Since $\kappa_1,\kappa_2\ne0$ and $\kappa_1+\kappa_2\neq0$, this proves \eqref{detMM} and completes the proof of (b).

Let $(u_1,v_1,u_2,v_2)\in\mathcal Y^k$ be as in \eqref{uvexp}. We now compute the blocks. For $j\ge2$, a direct computation gives \eqref{Mjinv}, so that
\begin{equation*}
a_j^i=\frac{A_j^i+B_j^i}{\kappa_i(j-1)},\qquad b_j^i=\frac{(2-j)A_j^i+jB_j^i}{\kappa_i(j-1)},\qquad j\ge2 .
\end{equation*}
For $j=1$, adding and subtracting the identities $A_1^i=\lambda_i+\kappa_ia_1^i$ and $B_1^i=\lambda_i-\kappa_ia_1^i$ we get
\begin{equation*}
\lambda_i=\frac{A_1^i+B_1^i}{2},\qquad a_1^i=\frac{A_1^i-B_1^i}{2\kappa_i},\qquad i=1,2 .
\end{equation*}
Finally, by \eqref{lambdai}, $(\beta_1,\beta_2)$ solves the linear system
\begin{equation}\label{deriv22}
\lambda_1=\beta_1x_0^*+\beta_2\Omega^*,\qquad \lambda_2=\beta_1(d-x_0^*)-\beta_2\Omega^*,
\end{equation}
whose determinant is $-d\Omega^*\neq0$. Adding its two equations gives $\beta_1=\frac{\lambda_1+\lambda_2}{d}$, and then the first equation gives
\begin{equation*}
\beta_2=\frac{\lambda_1-x_0^*\beta_1}{\Omega^*}=\frac{(d-x_0^*)\lambda_1-x_0^*\lambda_2}{d\,\Omega^*}.
\end{equation*}
This proves \eqref{MM1inv}. Using $x_0^*=\frac{d\kappa_2}{\kappa_1+\kappa_2}$, $d-x_0^*=\frac{d\kappa_1}{\kappa_1+\kappa_2}$ and $\Omega^*=\frac{\kappa_1+\kappa_2}{d^2}$, we have $\frac{d-x_0^*}{d\Omega^*}=\frac{d^2\kappa_1}{(\kappa_1+\kappa_2)^2}$ and $\frac{x_0^*}{d\Omega^*}=\frac{d^2\kappa_2}{(\kappa_1+\kappa_2)^2}$. Hence, the inverse of $D_\gfr\bfF(0,\gfr_0)$ is given explicitly by
\begin{equation}\label{g-1}
\begin{split}
\beta_1&\fit[.85\hsize]{{}=\frac{A_1^1+B_1^1+A_1^2+B_1^2}{2d},\qquad \beta_2=\frac{d^2\big(\kappa_1(A_1^1+B_1^1)-\kappa_2(A_1^2+B_1^2)\big)}{2(\kappa_1+\kappa_2)^2},}\\
h_i&=\frac{A_1^i-B_1^i}{2\kappa_i}\cos x+\sum_{j=2}^\infty\frac{A_j^i+B_j^i}{\kappa_i(j-1)}\cos(jx),\\
\eta_i&=-\frac{A_1^i-B_1^i}{2\kappa_i}\cos x+\sum_{j=2}^\infty\frac{(2-j)A_j^i+jB_j^i}{\kappa_i(j-1)}\cos(jx).
\end{split}
\end{equation}
Thus,  the inverse \eqref{g-1} can be written as  $D_\gfr\bfF(0,\gfr_0)^{-1}(u_1,v_1,u_2,v_2)=(\beta_1,\beta_2,h_1,\eta_1,h_2,\eta_2)$.

 Since $\frac{1}{j-1}\le\frac2j$ and $\frac{j}{j-1}\le2$ for $j\ge2$, \eqref{g-1} gives
\begin{equation*}
|a_j^i|\le\frac{2}{|\kappa_i|\,j}\big(|A_j^i|+|B_j^i|\big),\qquad |b_j^i|\le\frac{2}{|\kappa_i|}\big(|A_j^i|+|B_j^i|\big),\qquad j\ge2,
\end{equation*}
and $|\beta_1|+|\beta_2|+|a_1^1|+|a_1^2|\le C\sum_{i=1}^2\big(|A_1^i|+|B_1^i|\big)$. In other words, the first row of each block $\big(M_j^{(i)}\big)^{-1}$ gains a factor $\frac1j$, which corresponds to the extra derivative of $h_i\in X^{k+1}$, while the second row is bounded, which corresponds to $\eta_i\in X^k$. Then,  we obtain
\begin{equation*}
|\beta_1|+|\beta_2|+\sum_{i=1}^2\Big(\norm{h_i}_{X^{k+1}}+\norm{\eta_i}_{X^k}\Big)\le C\sum_{i=1}^2\Big(\norm{u_i}_{Y^k}+\norm{v_i}_{X^k}\Big),
\end{equation*}
and $(h_i,\eta_i)\in\mathcal A^k$ because $b_1^i=-a_1^i$. Hence $D_\gfr\bfF(0,\gfr_0)$ is a bijection from $\mathcal X^k$ onto $\mathcal Y^k$ with bounded inverse. Since it is bounded by Lemma \ref{lem:C1}, it is an isomorphism, which proves (c).
\end{proof}

\subsection{The angular velocity, the centroid and a priori estimates}
By Lemma \ref{critical}, the point $(0,\gfr_0)=(0,\Omega^*,x_0^*,0,0,0,0)$, with $\Omega^*$ and $x_0^*$ given in \eqref{vel}, is a trivial solution of $\bfF(\ep,\gfr)=0$, which corresponds to the co-rotating pair of point vortices $2\pi\kappa_1\boldsymbol\delta_0$ and $2\pi\kappa_2\boldsymbol\delta_{d\bfe_1}$. The angular velocity and the centroid are determined by the first Fourier coefficients of the equations. For $u\in Y^k$, $v\in X^k$ and $j\in\N^*$ we denote by $\mathcal R_i:=\mathcal R_{i1}+\mathcal R_{i2}+\mathcal R_{i3}$ and $\mathcal Q_i:=\mathcal Q_{i1}+\mathcal Q_{i2}+\mathcal Q_{i3}$, where $\mathcal R_{il}$ and $\mathcal Q_{il}$ are the remainders terms of Lemma \ref{lem:cont}. 

By \eqref{Fexp}--\eqref{Hexp}, for $i=1,2$, we have
\begin{equation}\label{FHrem}
\begin{split}
\mathcal F_i&=\Big(\Omega x_{0,i}-\frac{\kappa_{3-i}}{d}\Big)\sin x-\frac{\kappa_i}{2}r_i'(x)\\
&\qquad-\kappa_i\fint\frac{s_i(y)\sin(x-y)}{4\sin^2(\frac{x-y}{2})}\dy+\ep\mathcal R_i,\\
\mathcal H_i&=\Big(\Omega x_{0,i}-\frac{\kappa_{3-i}}{d}\Big)\cos x-\kappa_ir_i(x)\\
&\qquad+\kappa_i\fint\frac{r_i(x)-r_i(y)}{4\sin^2(\frac{x-y}{2})}\dy+\frac{\kappa_i}{2}s_i(x)+\ep\mathcal Q_i .
\end{split}
\end{equation}

To derive co-rotating asymmetric vortex sheet solutions using the implicit function theorem, we must select an appropriate $\Omega=\Omega(\ep,\Omega,x_0,r_1,s_1,r_2,s_2)$ that ensures the operator $\mathcal F(\ep,\Omega,x_0,r_1,s_1,r_2,s_2)$ maps $B_\varrho\subset\mathcal X^k$ into $\mathcal Y^k$. The following lemma establishes this choice

\begin{lemma}\label{lem2-5}
Let $(\ep,\gfr)\in(-\ep_1,\ep_1)\times\R\times\R\times B_\varrho$ with $r_i=\displaystyle\sum_{j\ge1}a_j^i\cos(jx)$ and $s_i=-a_1^i\cos x+\displaystyle\sum_{j\ge2}b_j^i\cos(jx)$, $i=1,2$. If $\bfF(\ep,\gfr)=0$, then, for $i=1,2$, we have
\begin{equation}\label{Omegaxb}
\Omega=\Omega^*+\ep\,\Omega_\ep(\ep,\gfr),\qquad \Omega_\ep:=-\frac{\mathfrak r_1+\mathfrak r_2}{d},
\end{equation}
and, if $\Omega\neq0$,
\begin{equation}\label{xbeps}
\fit[\hsize]{x_0=x_0^*+\ep\,x_{0,\ep}(\ep,\gfr),\qquad x_{0,\ep}:=-\frac{\frac{\kappa_2}{d}\Omega_\ep+\Omega^*\mathfrak r_1}{\Omega^*\,\Omega}.}
\end{equation}
where
\begin{equation*}
\fit[\hsize]{\mathfrak r_i:=\fint\big(\mathcal R_i(x)\sin x+\mathcal Q_i(x)\cos x\big)\dx,\qquad \mathfrak q_i:=\fint\big(\mathcal R_i(x)\sin x-\mathcal Q_i(x)\cos x\big)\dx .}
\end{equation*}
then $\mathcal{F^*}(\varepsilon,x_0,r_1,s_1,r_2,s_2):=\mathcal F(\varepsilon,\Omega(\varepsilon,g),x_0,r_1,s_1,r_2,s_2)$ maps $ (-\ep_1,\ep_1)\times\R\times\R\times B_\varrho$ into $\mathcal Y^k$.
It also holds that $\Omega_\varepsilon(\varepsilon,g)$ and $x_{0,\ep}(\varepsilon,g)$ remain continuously differentiable for all $(\varepsilon,g)\in (-\ep_1,\ep_1)\times\R\times\R\times B_\varrho$.
\end{lemma}

\begin{proof}
We first prove (a). Since $(r_i,s_i)\in\mathcal A^k$, the first Fourier coefficient of $s_i$ is $-a_1^i$, so that $r_i$ and $s_i$ can be written as in the statement. By \eqref{HLcos},
\begin{equation*}
\begin{gathered}
\fit[\hsize]{r_i'(x)=-\sum_{j\ge1}ja_j^i\sin(jx),\qquad \fint\frac{s_i(y)\sin(x-y)}{4\sin^2(\frac{x-y}{2})}\dy=-\frac{a_1^i}{2}\sin x+\frac12\sum_{j\ge2}b_j^i\sin(jx),}\\
\fint\frac{r_i(x)-r_i(y)}{4\sin^2(\frac{x-y}{2})}\dy=\frac12\sum_{j\ge1}ja_j^i\cos(jx),
\end{gathered}
\end{equation*}
and therefore
\begin{equation}\label{lincoef}
\begin{split}
-\frac{\kappa_i}{2}r_i'(x)-\kappa_i\fint\frac{s_i(y)\sin(x-y)}{4\sin^2(\frac{x-y}{2})}\dy&\fit[.85\hsize]{{}=\kappa_ia_1^i\sin x+\frac{\kappa_i}{2}\sum_{j\ge2}\big(ja_j^i-b_j^i\big)\sin(jx),}\\
-\kappa_ir_i(x)+\kappa_i\fint\frac{r_i(x)-r_i(y)}{4\sin^2(\frac{x-y}{2})}\dy+\frac{\kappa_i}{2}s_i(x)&\fit[.85\hsize]{{}=-\kappa_ia_1^i\cos x+\frac{\kappa_i}{2}\sum_{j\ge2}\big((j-2)a_j^i+b_j^i\big)\cos(jx).}
\end{split}
\end{equation}
Multiplying the first identity in \eqref{FHrem} by $\sin x$, the second one by $\cos x$, integrating over $[0,2\pi]$ and using \eqref{lincoef} and $\fint\sin^2x\dx=\fint\cos^2x\dx=\frac12$, we obtain
\begin{equation*}
\begin{split}
\fint\mathcal F_i(x)\sin x\dx&\fit[.85\hsize]{{}=\frac12\Big(\Omega x_{0,i}-\frac{\kappa_{3-i}}{d}\Big)+\frac{\kappa_i}{2}a_1^i+\ep\fint\mathcal R_i(x)\sin x\dx,}\\
\fint\mathcal H_i(x)\cos x\dx&\fit[.85\hsize]{{}=\frac12\Big(\Omega x_{0,i}-\frac{\kappa_{3-i}}{d}\Big)-\frac{\kappa_i}{2}a_1^i+\ep\fint\mathcal Q_i(x)\cos x\dx .}
\end{split}
\end{equation*}
If $\bfF(\ep,\gfr)=0$, both left-hand sides vanish. Adding the two identities we get $\Omega x_{0,i}-\frac{\kappa_{3-i}}{d}=-\ep\,\mathfrak r_i$, and subtracting them we get $\kappa_ia_1^i=-\ep\,\mathfrak q_i$. Adding now the first identity for $i=1$ and $i=2$, and using $x_{0,1}+x_{0,2}=x_0+(d-x_0)=d$, we obtain
\begin{equation*}
\Omega d=\frac{\kappa_1+\kappa_2}{d}-\ep\big(\mathfrak r_1+\mathfrak r_2\big),
\end{equation*}
that is, by \eqref{vel}, $\Omega=\Omega^*-\frac{\ep}{d}(\mathfrak r_1+\mathfrak r_2)$, which is \eqref{Omegaxb}. If $\Omega\ne0$, the first identity  with $i=1$ gives $x_0=\frac1\Omega\big(\frac{\kappa_2}{d}-\ep\,\mathfrak r_1\big)$, and, since $x_0^*=\frac{\kappa_2}{d\,\Omega^*}$ by \eqref{vel}, we have
\begin{equation*}
\fit[\hsize]{x_0-x_0^*=\frac{\frac{\kappa_2}{d}\big(\Omega^*-\Omega\big)-\ep\,\Omega^*\mathfrak r_1}{\Omega^*\,\Omega}=-\ep\,\frac{\frac{\kappa_2}{d}\Omega_\ep+\Omega^*\mathfrak r_1}{\Omega^*\,\Omega},}
\end{equation*}
which is \eqref{xbeps}. Finally,  as discussed in Lemma \ref{lem:C1}, it easily follows that the operators $\Omega_\varepsilon(\varepsilon,g)$ and $x_{0,\ep}(\varepsilon,g)$ are of class $C^1$ for all $(\varepsilon,g)\in (-\ep_1,\ep_1)\times\R\times\R\times B_\varrho$. This completes the proof.
\end{proof}

%%%
%%%
We shall also need the expansion, up to order $\ep^2$, of the interaction terms when the sheets are circles with uniform strength. To do this, we apply the Taylor formula \eqref{taylor}.

\begin{lemma}\label{lem:circles}
Let $|\ep|<\ep_1$ and $i\in\{1,2\}$, and assume that $r_1=s_1=r_2=s_2=0$. Then
\begin{equation}\label{circles}
\begin{gathered}
\mathcal F_{i3}=-\frac{\kappa_{3-i}}{d}\Big(\sin x+\frac{\ep b_i}{d}\sin(2x)\Big)+\ep^2\mathcal R_{i3}^{c},\\
\mathcal H_{i3}=-\frac{\kappa_{3-i}}{d}\Big(\cos x+\frac{\ep b_i}{d}\cos(2x)\Big)+\ep^2\mathcal Q_{i3}^{c},
\end{gathered}
\end{equation}
where $\mathcal R_{i3}^{c}\in Y^k$ and $\mathcal Q_{i3}^{c}\in X^k$ are given in by
\begin{equation}\label{RQc}
\begin{split}
\mathcal R_{i3}^{c}:=&\fit[.85\hsize]{\ \kappa_{3-i}\Big(\frac{b_i^2+b_{3-i}^2}{d^3}\sin x-d\sin x\fint\mathfrak K_i(x,y)\dy\Big),}\\
\tilde{\mathcal Q}_{i3}^{c}:=&\fit[.85\hsize]{\ \kappa_{3-i}\Big(\frac{3b_i^2+2b_{3-i}^2}{d^3}\cos x-d\cos x\fint\mathfrak K_i(x,y)\dy+\ep b_i\Big(\fint\mathfrak K_i(x,y)\dy-\frac{b_i^2+b_{3-i}^2}{d^4}\Big)}\\
&\fit[.85\hsize]{\qquad+\ep b_{3-i}\Big(\fint\cos(x-y)\mathfrak K_i(x,y)\dy-\frac{b_ib_{3-i}}{d^4}\Big)\Big),\qquad \mathcal Q_{i3}^{c}:=(I-\Pm)\tilde{\mathcal Q}_{i3}^{c}.}
\end{split}
\end{equation}
where they also satisfy $\norm{\mathcal Q_{i3}^{c}}_{X^k}+\norm{\mathcal R_{i3}^{c}}_{Y^k}\le C$.
\end{lemma}

\begin{proof}
Since $r_1=s_1=r_2=s_2=0$, we have $1+\ep b_jr_j=1+\ep b_js_j=1$ for $j=1,2$ and $r_i'=0$, so that, by \eqref{F4} and \eqref{H3}, we have
\begin{equation}\label{circ1}
\begin{split}
\mathcal F_{i4}=&\fit[.85\hsize]{{}-\kappa_{3-i}\,d\sin x\fint\frac{\dy}{\big(\ep b_i\cos x+\ep b_{3-i}\cos y-d\big)^2+\big(\ep b_i\sin x+\ep b_{3-i}\sin y\big)^2}}\\
\tilde{\mathcal H}_{i3}=&\fit[.85\hsize]{{}-\kappa_{3-i}\,d\cos x\fint\frac{\dy}{\big(\ep b_i\cos x+\ep b_{3-i}\cos y-d\big)^2+\big(\ep b_i\sin x+\ep b_{3-i}\sin y\big)^2}}\\
&\fit[.85\hsize]{{}+\kappa_{3-i}\,\ep b_i\fint\frac{\dy}{\big(\ep b_i\cos x+\ep b_{3-i}\cos y-d\big)^2+\big(\ep b_i\sin x+\ep b_{3-i}\sin y\big)^2}}\\
&\fit[.85\hsize]{{}+\kappa_{3-i}\,\ep b_{3-i}\fint\frac{\cos(x-y)}{\big(\ep b_i\cos x+\ep b_{3-i}\cos y-d\big)^2+\big(\ep b_i\sin x+\ep b_{3-i}\sin y\big)^2}\dy .}
\end{split}
\end{equation}
Expanding the squares in the denominator and using $\cos^2+\sin^2=1$ and $\cos x\cos y+\sin x\sin y=\cos(x-y)$, we get
\begin{equation}\label{circB}
\big(\ep b_i\cos x+\ep b_{3-i}\cos y-d\big)^2+\big(\ep b_i\sin x+\ep b_{3-i}\sin y\big)^2=d^2+B_i(x,y),
\end{equation}
where
\begin{equation*}
B_i(x,y):=-2d\ep\big(b_i\cos x+b_{3-i}\cos y\big)+\ep^2\big(b_i^2+b_{3-i}^2+2b_ib_{3-i}\cos(x-y)\big).
\end{equation*}
We now apply the Taylor formula \eqref{taylor} twice. First, with $\lambda=1$, $A=d^2$ and $B=B_i(x,y)$, we get
\begin{equation*}
\frac{1}{d^2+B_i(x,y)}=\frac{1}{d^2}-\int_0^1\frac{B_i(x,y)}{\big(d^2+tB_i(x,y)\big)^2}\dt ;
\end{equation*}
then, for each $t\in[0,1]$, with $\lambda=2$, $A=d^2$ and $B=tB_i(x,y)$,
\begin{equation*}
\frac{1}{\big(d^2+tB_i(x,y)\big)^2}=\frac{1}{d^4}-2\int_0^1\frac{tB_i(x,y)}{\big(d^2+stB_i(x,y)\big)^3}\,\mathrm{d}s .
\end{equation*}
Inserting the second identity in the first one, we arrive at
\begin{equation}\label{circTaylor}
\frac{1}{d^2+B_i(x,y)}=\frac{1}{d^2}-\frac{B_i(x,y)}{d^4}+\ep^2\mathfrak K_i(x,y),
\end{equation}
where
\begin{equation*}
\fit[\hsize]{\mathfrak K_i(x,y):=2\Big(\frac{B_i(x,y)}{\ep}\Big)^2\int_0^1\!\!\int_0^1\frac{t}{\big(d^2+stB_i(x,y)\big)^3}\,\mathrm{d}s\dt .}
\end{equation*}
The function $\mathfrak K_i$ is not singular at $\ep=0$, because
\begin{equation*}
\frac{B_i(x,y)}{\ep}=-2d\big(b_i\cos x+b_{3-i}\cos y\big)+\ep\big(b_i^2+b_{3-i}^2+2b_ib_{3-i}\cos(x-y)\big).
\end{equation*}
The definition of $B_i$ gives
\begin{equation*}
\begin{gathered}
\fint B_i(x,y)\dy=-2d\ep b_i\cos x+\ep^2\big(b_i^2+b_{3-i}^2\big),\qquad \fint\sin(x-y)B_i(x,y)\dy=-d\ep b_{3-i}\sin x,\\
\fint\cos(x-y)B_i(x,y)\dy=-d\ep b_{3-i}\cos x+\ep^2b_ib_{3-i} .
\end{gathered}
\end{equation*}
Therefore, by \eqref{circB} and \eqref{circTaylor},
\begin{equation}\label{circmeans}
\begin{split}
\fint\frac{\dy}{d^2+B_i(x,y)}&\fit[.85\hsize]{{}=\frac{1}{d^2}+\frac{2\ep b_i}{d^3}\cos x-\frac{\ep^2\big(b_i^2+b_{3-i}^2\big)}{d^4}+\ep^2\fint\mathfrak K_i(x,y)\dy,}\\
\fint\frac{\cos(x-y)}{d^2+B_i(x,y)}\dy&\fit[.85\hsize]{{}=\frac{\ep b_{3-i}}{d^3}\cos x-\frac{\ep^2b_ib_{3-i}}{d^4}+\ep^2\fint\cos(x-y)\mathfrak K_i(x,y)\dy .}
\end{split}
\end{equation}
Inserting \eqref{circmeans} in \eqref{circ1} and using $2\sin x\cos x=\sin(2x)$ and $2\cos^2x-1=\cos(2x)$, we obtain
\begin{equation*}
\begin{split}
\mathcal F_{i3}&=-\frac{\kappa_{3-i}}{d}\sin x-\frac{2\kappa_{3-i}\ep b_i}{d^2}\sin x\cos x+\ep^2\mathcal R_{i3}^{c}\\
&=-\frac{\kappa_{3-i}}{d}\Big(\sin x+\frac{\ep b_i}{d}\sin(2x)\Big)+\ep^2\mathcal R_{i3}^{c},\\
\tilde{\mathcal H}_{i3}&\fit[.85\hsize]{{}=-\frac{\kappa_{3-i}}{d}\cos x-\frac{2\kappa_{3-i}\ep b_i}{d^2}\cos^2x+\frac{\kappa_{3-i}\ep b_i}{d^2}+\ep^2\tilde{\mathcal Q}_{i3}^{c}}\\
&=-\frac{\kappa_{3-i}}{d}\Big(\cos x+\frac{\ep b_i}{d}\cos(2x)\Big)+\ep^2\tilde{\mathcal Q}_{i3}^{c},
\end{split}
\end{equation*}
where
\begin{equation}
\begin{split}
\mathcal R_{i4}^{c}:=&\fit[.85\hsize]{\ \kappa_{3-i}\Big(\frac{b_i^2+b_{3-i}^2}{d^3}\sin x-d\sin x\fint\mathfrak K_i(x,y)\dy\Big),}\\
\tilde{\mathcal Q}_{i3}^{c}:=&\fit[.85\hsize]{\ \kappa_{3-i}\Big(\frac{3b_i^2+2b_{3-i}^2}{d^3}\cos x-d\cos x\fint\mathfrak K_i(x,y)\dy+\ep b_i\Big(\fint\mathfrak K_i(x,y)\dy-\frac{b_i^2+b_{3-i}^2}{d^4}\Big)}\\
&\fit[.85\hsize]{\qquad+\ep b_{3-i}\Big(\fint\cos(x-y)\mathfrak K_i(x,y)\dy-\frac{b_ib_{3-i}}{d^4}\Big)\Big),\qquad \mathcal Q_{i3}^{c}:=(I-\Pm)\tilde{\mathcal Q}_{i3}^{c}.}
\end{split}
\end{equation}
Applying $I-\Pm$ to $\tilde{\mathcal H}_{i3}$ we obtain the expansion of $\mathcal H_{i3}$ in \eqref{circles}. Moreover, $\mathcal F_{i3}$ is odd and $\tilde{\mathcal H}_{i3}$ is even, as shown in the proof of Lemma \ref{lem:cont}, and so are the explicit terms in \eqref{circles}. Therefore $\mathcal R_{i4}^{c}$ is odd and $\mathcal Q_{i3}^{c}$ is even with zero mean, that is, $\mathcal R_{i4}^{c}\in Y^k$, $\mathcal Q_{i3}^{c}\in X^k$, and $\norm{\mathcal R_{i4}^{c}}_{Y^k}+\norm{\mathcal Q_{i3}^{c}}_{X^k}\le C$.
\end{proof}

\subsection{Existence of co-rotating asymmetric vortex sheets}

Now we provide a detailed version of the first part of Theorem \ref{thm:informal}.

\begin{theorem}\label{main}
Let  $d>0$, $b_1,b_2>0$ and $\kappa_1,\kappa_2\in\R\setminus\{0\}$ with $\kappa_1+\kappa_2\neq0$. Then the following assertions hold true.
\begin{enumerate}
\item There exist $\ep_0>0$, a neighborhood $\mathcal U_0$ of $\gfr_0$ in $\R\times\R\times B_\varrho$ and a unique $C^1$ function
\begin{equation*}
\gfr=(\Omega,x_0,r_1,s_1,r_2,s_2):[-\ep_0,\ep_0]\longrightarrow\mathcal U_0
\end{equation*}
such that $\bfF\big(\ep,\gfr(\ep)\big)=0$. Moreover, every zero of $\bfF$ in $[-\ep_0,\ep_0]\times\mathcal U_0$ belongs to this curve. Furthermore, $\gfr(0)=(\Omega^*,x_0^*,0,0,0,0)$, that is, the curve passes through the point vortex pair.
\item For each $\ep\in[-\ep_0,\ep_0]\setminus\{0\}$ and each $i\in\{1,2\}$, it always holds
\begin{equation*}
\big(r_i(\ep),s_i(\ep)\big)\neq(0,0).
\end{equation*}
That is, neither of the two sheets is a circle with uniform strength.
\item For each $\ep\in[-\ep_0,\ep_0]$, we have $\Omega(-\ep)=\Omega(\ep)$, $x_0(-\ep)=x_0(\ep)$ and
\begin{equation*}
r_i(-\ep)(x)=-r_i(\ep)(x+\pi),\qquad s_i(-\ep)(x)=-s_i(\ep)(x+\pi),\qquad i=1,2 .
\end{equation*}
In particular, the curves and the vortex strengths given by $-\ep$ coincide with those given by $\ep$, up to the shift $x\mapsto x+\pi$ of the parameter.
\item For each $\ep\in[-\ep_0,\ep_0]\setminus\{0\}$, the curves $\Gamma_i^\ep$ defined by \eqref{curves} and \eqref{ansatz} are real-analytic and enclose strictly convex domains, and the vortex strengths $\gamma_i$ are real-analytic.
\end{enumerate}
Consequently, for $0<\ep\le\ep_0$, the vorticity $\vartheta_\ep(t,X)=\vartheta_{0,\ep}\big(Q_{-\Omega(\ep)t}(X-x_0(\ep)\bfe_1)+x_0(\ep)\bfe_1\big)$ is a co-rotating asymmetric vortex sheet pair solution of \eqref{1-1}, rotating about $(x_0(\ep),0)$ with angular velocity $\Omega(\ep)$.
\end{theorem}

\begin{proof}
We first prove (1). By Lemmas \ref{lem:cont} and \ref{lem:C1}, the functional $\bfF:(-\ep_1,\ep_1)\times\R\times\R\times B_\varrho\to\mathcal Y^k$ is a $C^1-$function. By Lemma \ref{critical}, $\bfF(0,\gfr_0)=0$, and by Lemma \ref{lem:iso}, $D_\gfr\bfF(0,\gfr_0):\mathcal X^k\to\mathcal Y^k$ is an isomorphism. Now we can apply the implicit function theorem  and obtain $\ep_0>0$, a neighborhood $\mathcal U_0$ of $\gfr_0$ and a unique $C^1$ parametrization $\gfr:[-\ep_0,\ep_0]\to\mathcal U_0$ of the zeros of $\bfF$ in $[-\ep_0,\ep_0]\times\mathcal U_0$, with $\gfr(0)=\gfr_0$. Without loss of generality, we may take $\mathcal U_0=\{|\Omega-\Omega^*|<\delta,\ |x_0-x_0^*|<\delta\}\times\{\norm{r_i}_{X^{k+1}}+\norm{s_i}_{X^k}<\delta,\ i=1,2\}$ for some small $\delta>0$. The estimate follows from Lemma \ref{lem2-5} (b), reducing $\ep_0$ if necessary.

Now we prove (2). Assume by contradiction that $(r_i(\ep),s_i(\ep))=(0,0)$ for some $i\in\{1,2\}$ and some $\ep\neq0$. Then, by Lemma \ref{lem:circles}, we have $\mathcal F_{i1}=\Omega(\ep)x_{0,i}(\ep)\sin x$ and $\mathcal F_{i2}=0$. On the other hand, $\mathcal F_{i3}$ depends on $(r_{3-i},s_{3-i})$ only through $1+\ep b_{3-i}r_{3-i}(y)$ and $1+\ep b_{3-i}s_{3-i}(y)$, so, by the mean value theorem and (1),
\begin{equation*}
\fit[\hsize]{\Big\Vert \mathcal F_{i3}\big(\ep,\gfr(\ep)\big)-\mathcal F_{i3}\big|_{r_1=s_1=r_2=s_2=0}\Big\Vert_{Y^k}\le C|\ep|\big(\norm{r_{3-i}(\ep)}_{X^{k+1}}+\norm{s_{3-i}(\ep)}_{X^k}\big)\le C\ep^2 .}
\end{equation*}
Therefore, by Lemma \ref{lem:circles}, and since the map $u\mapsto2\fint u(x)\sin(2x)\dx$ is bounded on $Y^k$, the coefficient of $\sin(2x)$ in $\mathcal F_i(\ep,\gfr(\ep))$ is
\begin{equation*}
-\frac{\kappa_{3-i}}{d}\,\frac{\ep b_i}{d}+O(\ep^2)=-\frac{\kappa_{3-i}b_i}{d^2}\,\ep+O(\ep^2),
\end{equation*}
which is different from zero for $0<|\ep|\le\ep_0$ if $\ep_0$ is chosen small enough, since $\kappa_{3-i}b_i\neq0$. This contradicts $\mathcal F_i(\ep,\gfr(\ep))=0$. Notice that, in particular, the circles with uniform strengths are never solutions for $\ep\ne0$. By Lemma \ref{lem:circles}, we have that $\bfF(\ep,\Omega,x_0,0,0,0,0)\neq0$ for all $(\Omega,x_0)$, because of the term with $\sin(2x)$ in \eqref{circles}, whose coefficient $-\frac{\kappa_{3-i}b_i}{d^2}\ep+O(\ep^2)$ does not vanish for $0<|\ep|$ small. This reflects the fact that each sheet is deformed by the strain generated by the other one.

Next, we prove (3). For a function $u$ we set $(Tu)(x):=-u(x+\pi)$ and $(\tau_\pi u)(x):=u(x+\pi)$. If $u=\sum_ja_j\cos(jx)$, then $Tu=\sum_j(-1)^{j+1}a_j\cos(jx)$. Hence $T$ is an isometry of $X^k$, and it preserves $\mathcal A^k$ because it does not change the first Fourier coefficients. We claim that
\begin{equation}\label{symm}
\bfF\big(-\ep,\Omega,x_0,Tr_1,Ts_1,Tr_2,Ts_2\big)=-\tau_\pi\bfF\big(\ep,\Omega,x_0,r_1,s_1,r_2,s_2\big).
\end{equation}
%%
%%%
%%%
We start with $\mathcal F_{i21}$.  By replacing $(\ep,r_i,s_i)$ by $(-\ep,Tr_i,Ts_i)$, where $(Tu)(x)=-u(x+\pi)$. Since $(Tr_i)'(x)=-r_i'(x+\pi)$, we obtain
\begin{equation*}
\begin{split}
&\mathcal F_{i21}(-\ep,Tr_i,Ts_i)(x)\\
&\fit[.85\hsize]{\quad=-\kappa_i\fint\frac{\big(1+\ep b_ir_i(x+\pi)\big)\big(-r_i'(x+\pi)\big)\big(1-\cos((x+\pi)-(y+\pi))\big)\big(1+\ep b_is_i(y+\pi)\big)}{\ep^2b_i^2\big(r_i(x+\pi)-r_i(y+\pi)\big)^2+4\big(1+\ep b_ir_i(x+\pi)\big)\big(1+\ep b_ir_i(y+\pi)\big)\sin^2\big(\frac{(x+\pi)-(y+\pi)}{2}\big)}\dy}\\
&\fit[.85\hsize]{\quad=\kappa_i\fint\frac{\big(1+\ep b_ir_i(x+\pi)\big)r_i'(x+\pi)\big(1-\cos((x+\pi)-(y+\pi))\big)\big(1+\ep b_is_i(y+\pi)\big)}{\ep^2b_i^2\big(r_i(x+\pi)-r_i(y+\pi)\big)^2+4\big(1+\ep b_ir_i(x+\pi)\big)\big(1+\ep b_ir_i(y+\pi)\big)\sin^2\big(\frac{(x+\pi)-(y+\pi)}{2}\big)}\dy .}
\end{split}
\end{equation*}
The last integrand is $2\pi$-periodic in $y$, so the change of variables $y'=y+\pi$, which maps $[0,2\pi]$ onto $[\pi,3\pi]$, does not change its mean value over a period. Therefore
\begin{equation*}
\begin{split}
&\mathcal F_{i21}(-\ep,Tr_i,Ts_i)(x)\\
&\fit[.85\hsize]{\quad=\kappa_i\fint\frac{\big(1+\ep b_ir_i(x+\pi)\big)r_i'(x+\pi)\big(1-\cos((x+\pi)-y')\big)\big(1+\ep b_is_i(y')\big)}{\ep^2b_i^2\big(r_i(x+\pi)-r_i(y')\big)^2+4\big(1+\ep b_ir_i(x+\pi)\big)\big(1+\ep b_ir_i(y')\big)\sin^2\big(\frac{(x+\pi)-y'}{2}\big)}\dy'}\\
&\quad=-\mathcal F_{i21}(\ep,r_i,s_i)(x+\pi).
\end{split}
\end{equation*}
In the same way we can obtain
\begin{equation}\label{F21}
\mathcal F_{i2}(-\ep,Tr_i,Ts_i)=-\tau_\pi\mathcal F_{i2}(\ep,r_i,s_i).
\end{equation}
%%%
%%%
Making the change of variables $y\to y+\pi$  and using $\frac{1}{(-\ep)b_i}=-\frac{1}{\ep b_i}$, we get
\begin{equation*}
\fit[\hsize]{\mathcal F_{i2}(-\ep,Tr_i,Ts_i)=-\tau_\pi \mathcal F_{i2}(\ep,r_i,s_i),\qquad \tilde{\mathcal H}_{i2}(-\ep,Tr_i,Ts_i)=-\tau_\pi\tilde{\mathcal H}_{i2}(\ep,r_i,s_i).}
\end{equation*}
For the term $\mathcal F_{i1}$ and  ${\mathcal H}_{i1}$, we use $\cos x=-\cos(x+\pi)$ and $\sin x=-\sin(x+\pi)$ to get
\begin{equation*}
\fit[\hsize]{\mathcal F_{i1}(-\ep,Tr_i,Ts_i)=-\tau_\pi \mathcal F_{i1}(\ep,r_i,s_i),\qquad \tilde{\mathcal H}_{i1}(-\ep,Tr_i,Ts_i)=-\tau_\pi\tilde{\mathcal H}_{i1}(\ep,r_i,s_i).}
\end{equation*}
 Finally,  we deal with  $\mathcal F_{i3}$, $\mathcal F_{i4}$  and $\tilde{\mathcal H}_{i3}$. Again by the change of  variables $y\to y+\pi$, we have
\begin{equation*}
\fit[\hsize]{\mathcal F_{i3}(-\ep,Tr_i,Ts_i)=-\tau_\pi \mathcal F_{i3}(\ep,r_i,s_i),\qquad \mathcal F_{i4}(-\ep,Tr_i,Ts_i)=-\tau_\pi \mathcal F_{i4}(\ep,r_i,s_i).}
\end{equation*}
and
\begin{equation*}
\fit[\hsize]{ \tilde{\mathcal H}_{i3}(-\ep,Tr_i,Ts_i)=-\tau_\pi\tilde{\mathcal H}_{i3}(\ep,r_i,s_i).}
\end{equation*}
Thus \eqref{symm} follows and we obtain (3).

Finally, we prove (4). Since $r_i(\ep)\in X^{k+1}$ and $s_i(\ep)\in X^k$, they are real-analytic, and so are the curves $\Gamma_i^\ep$ and the strengths $\gamma_i$. To prove the convexity, we compute the signed curvature of $\Gamma_i^\ep$. The formula for the curvature of a polar curve gives
\begin{equation*}
\begin{split}
|\ep|b_i\,\mathfrak c_i(x)&=\frac{R_i(x)^2+2R_i'(x)^2-R_i(x)R_i''(x)}{\big(R_i(x)^2+R_i'(x)^2\big)^{\frac32}}\\
&\fit[.85\hsize]{{}=\frac{(1+\ep b_ir_i(x))^2+2\ep^2b_i^2r_i'(x)^2-\ep b_ir_i''(x)(1+\ep b_ir_i(x))}{\big((1+\ep b_ir_i(x))^2+\ep^2b_i^2r_i'(x)^2\big)^{\frac32}}=\frac{1+O(\ep)}{1+O(\ep)}>0,}
\end{split}
\end{equation*}
for $\ep\neq0$ small and each $x\in[0,2\pi)$. Hence the signed curvature is strictly positive and each $\Gamma_i^\ep$ encloses a strictly convex domain.

Finally, by the derivation of \eqref{rot}--\eqref{rotT} in Section \ref{sec:rot} and \eqref{ansatz}, $\bfF(\ep,\gfr(\ep))=0$ is equivalent to \eqref{normal}--\eqref{tangential}. Then $\vartheta_\ep$ is a co-rotating vortex sheet solution. This completes the proof of Theorem \ref{main}.
\end{proof}

%%%
%%%%
%%%%

\section{Existence of  traveling asymmetric vortex sheets}\label{sec4}

In this section we look for traveling solutions of the form \eqref{travsol}. By \eqref{t3}--\eqref{t3T} and \eqref{ansatz}, we have to find zeros of
\begin{equation*}
\bfG(\ep,U,\kappa_2,r_1,s_1,r_2,s_2):=\big(\mathcal G_1,\mathcal J_1,\mathcal G_2,\mathcal J_2\big),
\end{equation*}
where, for $i=1,2$,
\begin{equation*}
\fit[\hsize]{\mathcal G_i:=\mathcal G_{i1}+\mathcal F_{i2}-\mathcal F_{i3}-\mathcal F_{i4},\qquad \mathcal J_i:=\mathcal J_{i1}+\mathcal H_{i2}-\mathcal H_{i3},\qquad \mathcal J_{i1}:=(I-\Pm)\tilde{\mathcal J}_{i1},}
\end{equation*}
with
\begin{align}
\mathcal G_{i1}&:=-U\Big(\big(1+\ep b_ir_i(x)\big)\sin x-\ep b_ir_i'(x)\cos x\Big),\label{G1}\\
\tilde{\mathcal J}_{i1}&\fit[.85\hsize]{:=-U\,\frac{1+\ep b_is_i(x)}{\Ufull{i}}\Big(\ep b_ir_i'(x)\sin x+\big(1+\ep b_ir_i(x)\big)\cos x\Big),}\label{J1}
\end{align}
and  $\mathcal F_{i2}, \mathcal F_{i3}, \mathcal F_{i4}, \mathcal H_{i2}$, and $\mathcal H_{i3}$ are the functionals defined in \eqref{F2}--\eqref{H3}. otice that the functional $\bfG$ is essentially identical to the functional $\bfF$ studied in the previous section, with only three notable modifications. First, the rotation terms $\mathcal F_{i1}$ and $\mathcal H_{i1}$ are replaced by the translation terms $\mathcal G_{i1}$ and $\mathcal J_{i1}$. Second, the interaction terms undergo a sign change since the two sheets possess circulations of opposite signs. Finally, the unknown parameters are now the speed $U$ and the circulation $\kappa_2$, rather than $\Omega$ and $x_0$. In particular, $\kappa_2$ enters linearly in $\mathcal F_{22}$, $\mathcal H_{22}$, $\mathcal F_{13}$ and $\mathcal H_{13}$. We set
\begin{equation*}
\gfr:=(U,\kappa_2,r_1,s_1,r_2,s_2),\qquad \gfr_0:=(U^*,\kappa_1,0,0,0,0),\qquad U^*=\frac{\kappa_1}{d},
\end{equation*}
and we use the spaces $\mathcal X^k$, $\mathcal Y^k$ and the neighborhood $B_\varrho$ of Section \ref{sec:exist}.

\begin{lemma}\label{lem:trav}
Let $\kappa_1\in\R\setminus\{0\}$. The following assertions hold true.
\begin{enumerate}
\item The functional $\bfG$ is of class $C^1$ from $(-\ep_1,\ep_1)\times\R\times\R\times B_\varrho$ into $\mathcal Y^k$, and for $i=1,2$,
\begin{equation}\label{Gexp}
\begin{split}
\mathcal G_i&=\Bigl(-U+\frac{\kappa_{3-i}}{d}\Bigr)\sin x-\frac{\kappa_i}{2}r_i'(x)-\kappa_i\fint\frac{s_i(y)\sin(x-y)}{4\sin^2(\frac{x-y}{2})}\dy\\
&\quad+\ep\big(\mathcal R_{i1}^{\rm tr}+\mathcal R_{i2}-\mathcal R_{i3}\big),\\
\mathcal J_i&=\Bigl(-U+\frac{\kappa_{3-i}}{d}\Bigr)\cos x-\kappa_ir_i(x)+\kappa_i\fint\frac{r_i(x)-r_i(y)}{4\sin^2(\frac{x-y}{2})}\dy\\
&\quad+\frac{\kappa_i}{2}s_i(x)+\ep\big(\mathcal Q_{i1}^{\rm tr}+\mathcal Q_{i2}-\mathcal Q_{i3}\big),
\end{split}
\end{equation}
where $\mathcal R_{i2},\mathcal R_{i3},\mathcal Q_{i2},\mathcal Q_{i3}$ are given in Lemma \ref{lem:cont} and
\begin{equation*}
\begin{gathered}
\mathcal R_{i1}^{\rm tr}:=-Ub_i\big(r_i(x)\sin x-r_i'(x)\cos x\big),\\
\fit[\hsize]{\mathcal Q_{i1}^{\rm tr}:=-Ub_i(I-\Pm)\left[\frac{\big(1+\ep b_is_i(x)\big)r_i'(x)\sin x+\big(s_i(x)-r_i(x)+\ep b_i(r_i(x)s_i(x)-r_i(x)^2-r_i'(x)^2)\big)\cos x}{\Ufull{i}}\right].}
\end{gathered}
\end{equation*}
\item We have $\bfG(0,U,\kappa_2,0,0,0,0)=0$ if and only if $\kappa_2=\kappa_1$ and $U=U^*$. In other words, the point vortices $2\pi\kappa_1\boldsymbol\delta_0$ and $-2\pi\kappa_1\boldsymbol\delta_{d\bfe_1}$ travel in the direction $\bfe_2$ with speed $U^*=\frac{\kappa_1}{d}$.
\item For all $\hfr=(\beta_1,\beta_2,h_1,\eta_1,h_2,\eta_2)\in\mathcal X^k$ with $h_i,\eta_i$ as in \eqref{hexp}, one has
\begin{equation}\label{DG}
\begin{split}
D_\gfr \mathcal G_i(0,\gfr_0)\hfr&\fit[.85\hsize]{{}=\Big(-\beta_1+\delta_{i1}\frac{\beta_2}{d}\Big)\sin x+\sum_{j=1}^\infty\hat A_j^i\sin(jx),}\\
D_\gfr \mathcal J_i(0,\gfr_0)\hfr&=\Big(-\beta_1+\delta_{i1}\frac{\beta_2}{d}\Big)\cos x+\sum_{j=1}^\infty\hat B_j^i\cos(jx),
\end{split}
\end{equation}
where $\delta_{11}=1$, $\delta_{21}=0$, and
\begin{equation*}
\begin{pmatrix}\hat A_j^i\\ \hat B_j^i\end{pmatrix}=M_j\begin{pmatrix}a_j^i\\ b_j^i\end{pmatrix},\qquad M_j:=\frac{\kappa_1}{2}\begin{pmatrix}j&-1\\ j-2&1\end{pmatrix}.
\end{equation*}
\item The linearized operator $D_\gfr\bfG(0,\gfr_0):\mathcal X^k\to\mathcal Y^k$ is an isomorphism.
\end{enumerate}
\end{lemma}

\begin{proof}
The proof of this lemma follows the same steps as in the co-rotating case, so we only point out the differences.

We first prove (1). The functionals $\mathcal F_{i2},\mathcal F_{i3},\mathcal H_{i2},\mathcal H_{i3}$ have been studied in Lemmas \ref{lem:cont} and \ref{lem:C1}, and they depend linearly on $\kappa_1,\kappa_2$. The terms $\mathcal G_{i1}$ and $\mathcal J_{i1}$ have the same structure as $\mathcal F_{i1}$ and $\mathcal H_{i1}$. In particular $\mathcal G_{i1}$ is odd, $\tilde{\mathcal J}_{i1}$ is even. Hence the functional $\bfG$ is of class $C^1$. The expansions \eqref{Gexp} follow from \eqref{Fexp}--\eqref{Hexp} and from
\begin{equation*}
\begin{gathered}
\mathcal G_{i1}=-U\sin x-U\ep b_i\big(r_i(x)\sin x-r_i'(x)\cos x\big),\\
\fit[\hsize]{\tilde{\mathcal J}_{i1}=-U\cos x-U\ep b_i\,\frac{\big(1+\ep b_is_i(x)\big)r_i'(x)\sin x+\big(s_i(x)-r_i(x)+\ep b_i(r_i(x)s_i(x)-r_i(x)^2-r_i'(x)^2)\big)\cos x}{\Ufull{i}}.}
\end{gathered}
\end{equation*}

Concerning (2), by \eqref{Gexp} with $\ep=0$ and $r_i=s_i=0$,
\begin{equation*}
\begin{gathered}
\mathcal G_1=\Big(-U+\frac{\kappa_2}{d}\Big)\sin x,\qquad \mathcal J_1=\Big(-U+\frac{\kappa_2}{d}\Big)\cos x,\\
\mathcal G_2=\Big(-U+\frac{\kappa_1}{d}\Big)\sin x,\qquad \mathcal J_2=\Big(-U+\frac{\kappa_1}{d}\Big)\cos x .
\end{gathered}
\end{equation*}
Therefore $\bfG(0,U,\kappa_2,0,0,0,0)=0$ if and only if $U=\frac{\kappa_2}{d}=\frac{\kappa_1}{d}$.

To prove (3), differentiating \eqref{Gexp} at $(0,\gfr_0)$, we get $\partial_UG_i=-\sin x$, $\partial_UJ_i=-\cos x$, $\partial_{\kappa_2}\mathcal G_1=\frac1d\sin x$, $\partial_{\kappa_2}\mathcal J_1=\frac1d\cos x$, while, by \eqref{Gexp},
\begin{equation*}
\begin{gathered}
\partial_{\kappa_2}\mathcal G_2(0,\gfr_0)=-\frac12r_2'(x)-\fint\frac{s_2(y)\sin(x-y)}{4\sin^2(\frac{x-y}{2})}\dy=0,\\
\partial_{\kappa_2}\mathcal J_2(0,\gfr_0)=-r_2(x)+\fint\frac{r_2(x)-r_2(y)}{4\sin^2(\frac{x-y}{2})}\dy+\frac12s_2(x)=0,
\end{gathered}
\end{equation*}
since $r_2=s_2=0$ at $\gfr_0$. The derivatives with respect to $(r_i,s_i)$ are given by \eqref{lin0} with $\kappa_2=\kappa_1$, and the computation in the proof of Lemma \ref{lem:iso} gives \eqref{DG}.

Finally, to prove (4), let $w=(u_1,v_1,u_2,v_2)\in\mathcal Y^k$ with $u_i=\displaystyle\sum_jA_j^i\sin(jx)$ and $v_i=\displaystyle\sum_jB_j^i\cos(jx)$. For $j\ge2$, $M_j$ is invertible since $\kappa_1\neq0$, and $M_j^{-1}$ is given by \eqref{Mjinv} with $\kappa_i$ replaced by $\kappa_1$. For $j=1$, using $b_1^i=-a_1^i$, we must solve
\begin{equation*}
\fit[\hsize]{\lambda_i+\kappa_1a_1^i=A_1^i,\qquad \lambda_i-\kappa_1a_1^i=B_1^i,\qquad \text{where}\quad \lambda_1:=-\beta_1+\frac{\beta_2}{d},\quad \lambda_2:=-\beta_1 .}
\end{equation*}
Thus $a_1^i=\frac{A_1^i-B_1^i}{2\kappa_1}$ and $\lambda_i=\frac{A_1^i+B_1^i}{2}$. The linear map $(\beta_1,\beta_2)\mapsto(\lambda_1,\lambda_2)$ has determinant $\frac1d\neq0$, for any value of $\kappa_1$, and its inverse is $\beta_1=-\lambda_2$, $\beta_2=d(\lambda_1-\lambda_2)$. Hence
\begin{equation*}
\fit[\hsize]{D_\gfr\bfG(0,\gfr_0)^{-1}w=\Big(-\frac{A_1^2+B_1^2}{2},\ \frac d2\big(A_1^1+B_1^1-A_1^2-B_1^2\big),\ h_1,\eta_1,h_2,\eta_2\Big),}
\end{equation*}
where $h_i,\eta_i$ are given by the second line of \eqref{g-1} with $\kappa_i$ replaced by $\kappa_1$. As in the proof of Lemma \ref{lem:iso}, the inverse is bounded from $\mathcal Y^k$ into $\mathcal X^k$, which completes the proof.
\end{proof}

Now we state the precise version of the second part of Theorem \ref{thm:informal}.

\begin{theorem}\label{trans}
Let $d>0$, $b_1,b_2>0$ and $\kappa_1\in\R\setminus\{0\}$. Then the following assertions hold true.
\begin{enumerate}
\item There exist $\ep_0>0$, a neighborhood $\mathcal U_0$ of $\gfr_0$ in $\R\times\R\times B_\varrho$ and a unique $C^1$ function
\begin{equation*}
\gfr=(U,\kappa_2,r_1,s_1,r_2,s_2):[-\ep_0,\ep_0]\longrightarrow\mathcal U_0
\end{equation*}
such that $\bfG(\ep,\gfr(\ep))=0$. Moreover, every zero of $\bfG$ in $[-\ep_0,\ep_0]\times\mathcal U_0$ belongs to this curve. Furthermore, $\gfr(0)=\gfr_0=(U^*,\kappa_1,0,0,0,0)$ with $U^*=\frac{\kappa_1}{d}$, that is,
\begin{equation*}
U(0)=\frac{\kappa_1}{d},\qquad \kappa_2(0)=\kappa_1,\qquad r_i(0)=s_i(0)=0,\quad i=1,2,
\end{equation*}
so that the curve passes through the pair of point vortices $2\pi\kappa_1\boldsymbol\delta_0$ and $-2\pi\kappa_1\boldsymbol\delta_{d\bfe_1}$, which travels in the direction $\bfe_2$ with speed $U^*$.
\item For each $\ep\in[-\ep_0,\ep_0]\setminus\{0\}$ and $i\in\{1,2\}$, one has $\big(r_i(\ep),s_i(\ep)\big)\neq(0,0)$.
\item For each $\ep\in[-\ep_0,\ep_0]$, $U(-\ep)=U(\ep)$, $\kappa_2(-\ep)=\kappa_2(\ep)$, and $r_i(-\ep)(x)=-r_i(\ep)(x+\pi)$, $s_i(-\ep)(x)=-s_i(\ep)(x+\pi)$, $i=1,2$.
\item For each $\ep\in[-\ep_0,\ep_0]\setminus\{0\}$, the curves $\Gamma_i^\ep$ are real-analytic and enclose strictly convex domains, and the vortex strengths $\gamma_i$ are real-analytic.
\end{enumerate}
Consequently, for $0<\ep\le\ep_0$, the vorticity $\vartheta_\ep(t,X)=\vartheta_{0,\ep}(X-tU(\ep)\bfe_2)$ is a traveling asymmetric vortex sheets  of \eqref{1-1}.
\end{theorem}

\begin{proof}
The proof follows the same steps as the proof of Theorem \ref{main}, so we only sketch it.

%%%%
%%%
%%%
The proof follows the same steps as the proof of Theorem \ref{main}, so we only sketch it.

We first prove (1). By Lemma \ref{lem:trav}, $\bfG$ is of class $C^1$, $\bfG(0,\gfr_0)=0$ and $D_\gfr\bfG(0,\gfr_0)$ is an isomorphism from $\mathcal X^k$ onto $\mathcal Y^k$. The implicit function theorem gives the curve $\gfr(\ep)$, with $\mathcal U_0$ a product of balls centered at $\gfr_0$.

Since $\mathcal F_{i2},\mathcal F_{i3},\mathcal F_{i4},\mathcal H_{i2},\mathcal H_{i3}$ depend linearly on $\kappa_1,\kappa_2$, Lemma \ref{lem:cont} and the explicit form of $\mathcal R_{i1}^{\rm tr}$ and $\mathcal Q_{i1}^{\rm tr}$ show that the remainders
\begin{equation*}
\fit[\hsize]{\mathcal R_i^{\rm tr}:=\mathcal R_{i1}^{\rm tr}+\mathcal R_{i2}-\mathcal R_{i3},\qquad \mathcal Q_i^{\rm tr}:=\mathcal Q_{i1}^{\rm tr}+\mathcal Q_{i2}-\mathcal Q_{i3}}
\end{equation*}
satisfy $\norm{\mathcal R_i^{\rm tr}}_{Y^k}+\norm{\mathcal Q_i^{\rm tr}}_{X^k}\le C_1$ along the curve, $i=1,2$. We write $r_i=\displaystyle\sum_{j\ge1}a_j^i\cos(jx)$ and $s_i=-a_1^i\cos x+\displaystyle\sum_{j\ge2}b_j^i\cos(jx)$. Multiplying the identities \eqref{Gexp} by $\sin x$ and $\cos x$, respectively, integrating over $[0,2\pi]$ and using \eqref{lincoef}, we obtain, as in the proof of Lemma \ref{lem2-5} (a),
\begin{equation*}
-U+\frac{\kappa_{3-i}}{d}=-\ep\,\mathfrak r_i^{\rm tr},\qquad \kappa_ia_1^i=-\ep\,\mathfrak q_i^{\rm tr},\qquad i=1,2,
\end{equation*}
where
\begin{equation*}
\fit[\hsize]{\mathfrak r_i^{\rm tr}:=\fint\big(\mathcal R_i^{\rm tr}(x)\sin x+\mathcal Q_i^{\rm tr}(x)\cos x\big)\dx,\qquad \mathfrak q_i^{\rm tr}:=\fint\big(\mathcal R_i^{\rm tr}(x)\sin x-\mathcal Q_i^{\rm tr}(x)\cos x\big).\dx}
\end{equation*}
Taking $i=2$ and then $i=1$, we find
\begin{equation*}
\fit[\hsize]{U=\frac{\kappa_1}{d}+\ep\,\mathfrak r_2^{\rm tr},\qquad \kappa_2=\kappa_1+\ep\,d\big(\mathfrak r_2^{\rm tr}-\mathfrak r_1^{\rm tr}\big),}
\end{equation*}
 which completes the proof of (1).

Let us now prove (2). Assume that $(r_i(\ep),s_i(\ep))=(0,0)$ for some $i$ and some $\ep\neq0$. Then $\mathcal G_{i1}=-U\sin x$ and $\mathcal F_{i2}=0$, and, as in the proof of Theorem \ref{main} (2), the functional $\mathcal F_{i4}(\ep,\gfr(\ep))$ differs from the expression \eqref{circles} by $O(\ep^2)$ in $Y^k$. Hence the coefficient of $\sin(2x)$ in $\mathcal G_i(\ep,\gfr(\ep))$ equals $\frac{\kappa_{3-i}b_i}{d^2}\ep+O(\ep^2)$, which does not vanish for $0<|\ep|\le\ep_0$ small. This is a contradiction.

Concerning (3), as in the proof of item (3) of Theorem \ref{main}, one checks that
\begin{equation*}
\bfG(-\ep,U,\kappa_2,Tr_1,Ts_1,Tr_2,Ts_2)=-\tau_\pi\bfG(\ep,U,\kappa_2,r_1,s_1,r_2,s_2) .
\end{equation*}
 The conclusion follows from the uniqueness in (1).

Moreover, item (4) follows by the exact same argument as the corresponding part of Theorem \ref{main}. 

Finally, by \eqref{ansatz} and the derivation of \eqref{t3}--\eqref{t3T} in Section \ref{sec:trav}, the relation $\bfG(\ep,\gfr(\ep))=0$ is equivalent to \eqref{travNT}, so $\vartheta_\ep$ is a traveling vortex sheet solution of \eqref{1-1}.
\end{proof}

\appendix

\section{The Hilbert transform on the torus}\label{app:hilbert}
In this appendix we collect the identities for the singular integrals which are used throughout the paper. All the integrals are understood in the principal value sense. These identities are classical consequences of the properties of the conjugate function (Hilbert transform) on the torus $\T=\R/2\pi\Z$. The first identity is 
\begin{equation}\label{HL}
\fint\frac{\sin(x-y)}{4\sin^2(\frac{x-y}{2})}\dy=0 .
\end{equation} 
Now, for every $j\in\N^*$, we have
\begin{equation}\label{HLcos}
\begin{gathered}
\fit[\hsize]{\fint\frac{\cos(jy)\sin(x-y)}{4\sin^2(\frac{x-y}{2})}\dy=\frac12\sin(jx),\qquad \fint\frac{\sin(jy)\sin(x-y)}{4\sin^2(\frac{x-y}{2})}\dy=-\frac12\cos(jx),}\\
\fit[\hsize]{\fint\frac{\cos(jx)-\cos(jy)}{4\sin^2(\frac{x-y}{2})}\dy=\frac j2\cos(jx),\qquad \fint\frac{\sin(jx)-\sin(jy)}{4\sin^2(\frac{x-y}{2})}\dy=\frac j2\sin(jx),}
\end{gathered}
\end{equation}
see \cite[proof of Lemma 3.3]{CQZ}. Consequently, for $u(x)=c_0+\displaystyle\sum_{j\ge1}\big(c_j\cos(jx)+d_j\sin(jx)\big)$ in $L^2(\T)$, we have
\begin{equation}\label{ident}
\begin{split}
\fint\frac{u(y)\sin(x-y)}{4\sin^2(\frac{x-y}{2})}\dy&=\frac12\sum_{j\ge1}\big(c_j\sin(jx)-d_j\cos(jx)\big),\\
\fint\frac{u(x)-u(y)}{4\sin^2(\frac{x-y}{2})}\dy&=\frac12\sum_{j\ge1}j\big(c_j\cos(jx)+d_j\sin(jx)\big) .
\end{split}
\end{equation}

%%%
%%%
%%

\section{The implicit function theorem}\label{app:ift}

In this appendix we recall the implicit function theorem in the form used in the proofs of Theorems \ref{main} and \ref{trans}. For the proof we refer to \cite[Theorem I.1.1]{Kie}.

\begin{theorem}[Implicit function theorem]\label{IFT}
Let $\mathcal X,\mathcal Y$ be Banach spaces, $I\subset\R$ an open interval containing $0$, $\mathcal U\subset\mathcal X$ an open set containing $\gfr_0$, and let $\bfF:I\times\mathcal U\to\mathcal Y$ be continuous with $\bfF(0,\gfr_0)=0$. Assume that the partial Fr\'echet derivative $D_\gfr\bfF$ exists and is continuous on $I\times\mathcal U$, and that $D_\gfr\bfF(0,\gfr_0):\mathcal X\to\mathcal Y$ is an isomorphism. Then there exist $\ep_0>0$ and a neighborhood $\mathcal U_0\subset\mathcal U$ of $\gfr_0$ such that the set of zeros of $\bfF$ in $[-\ep_0,\ep_0]\times\mathcal U_0$ is a curve $\{(\ep,\gfr(\ep)):\ |\ep|\le\ep_0\}$, with $\gfr:[-\ep_0,\ep_0]\to\mathcal U_0$ continuous and $\gfr(0)=\gfr_0$. If, moreover, $\bfF\in C^1(I\times\mathcal U;\mathcal Y)$, then $\gfr\in C^1$, and $\gfr'(\ep)=-D_\gfr\bfF(\ep,\gfr(\ep))^{-1}\partial_\ep\bfF(\ep,\gfr(\ep))$.
\end{theorem}

\end{document}